\documentclass[a4paper,fleqn]{cas-sc}

\usepackage[authoryear]{natbib}
\usepackage{amsmath,amssymb}
\usepackage{algorithm}
\usepackage{algpseudocode}
\usepackage{booktabs}
\usepackage{longtable}
\usepackage[version=3]{mhchem}
\usepackage{tablefootnote}
\usepackage{array}
\usepackage{siunitx}
\usepackage{lineno}

\def\tsc#1{\csdef{#1}{\textsc{\lowercase{#1}}\xspace}}
\tsc{WGM}
\tsc{QE}

\newcommand{\Tair}{\mathrm{T,air}}
\newcommand{\COair}{\mathrm{CO_{2},air}}

\newcommand{\piRL}{\pi_{\mathrm{RL}}}

\begin{document}
\let\WriteBookmarks\relax
\def\floatpagepagefraction{1}
\def\textpagefraction{.001}

\shorttitle{Trajectory-selection RL--MPC for greenhouse fruit-production}    

\shortauthors{Van Laatum B et al.}

\title [mode = title]{Improving greenhouse fruit-production control by integrating reinforcement learning into short-horizon model predictive control}



%

\author[1]{B. van Laatum}[type=editor, 
    orcid=0009-0003-4675-836X]

\cormark[1]
\ead{bart.vanlaatum@wur.nl}



\credit{Conceptualization, Methodology, Software, Formal Analysis, Writing - Original Draft, Writing - Review \& Editing, Visualization, Project administration}

\affiliation[1]{organization={Agricultural Biosystems Engineering, Wageningen University \& Research},
                addressline={Droevendaalsesteeg 4}, 
                city={Wageningen},
                postcode={6708 PB}, 
                country={The Netherlands}}

\author[2]{S. Msaad}[type=editor,
orcid=0009-0002-2883-6048]

\credit{Conceptualization, Methodology, Writing - Review \& Editing}

\affiliation[2]{organization={Delft Center for Systems and Control (DCSC), Delft University of Technology},
                addressline={Mekelweg 5}, 
                city={Delft},
                postcode={2628 CD}, 
                country={The Netherlands}}

\author[1]{E. J. van Henten}[type=editor, 
    orcid=0000-0002-1623-9855]
\credit{Supervision, Writing - Review \& Editing}

\author[2]{R.D. McAllister}[type=editor, 
    orcid=0000-0002-5687-6875]
\credit{Conceptualization, Methodology, Writing - Review \& Editing, Supervision}

\author[1, 3]{S. Boersma}[type=editor, 
    orcid=0000-0002-3117-6514]
\credit{Conceptualization, Methodology, Writing - Review \& Editing, Supervision}

\affiliation[3]{organization={Biometris, Wageningen Research},
                addressline={Droevendaalsesteeg 4}, 
                city={Wageningen},
                postcode={6708 PB}, 
                country={The Netherlands}}

\cortext[1]{Corresponding author}


\begin{abstract}
Greenhouse fruit-production control aims to maximize the economic performance (fruit revenue minus operating costs) while operating within system constraints under external weather disturbances.
Control methods need to balance the delayed economic benefit of fruit yield with current operating costs.
For such problems, model predictive control (MPC) can explicitly handle system constraints under future weather disturbances, but can become computationally demanding when using sufficiently long prediction horizons for (relatively large) nonlinear greenhouse fruit production models.
In contrast, reinforcement learning (RL) can learn control policies offline while considering longer-term economic performance, but struggles to enforce system constraints, and performance may degrade under unseen weather trajectories.
This work proposes \textit{trajectory-selection RL--MPC}, a framework that incorporates longer-term economic information of fruit yield into a short-horizon MPC optimization problem.
The framework uses an RL rollout trajectory to define a terminal region constraint and terminal cost.
Next, a nonlinear MPC solves a short-horizon optimization problem with these terminal ingredients to find a local optimum.
Finally, the framework selects and executes the first input from the trajectory with the better objective value, either from the MPC-predicted or the RL rollout trajectory.
The method is applied to GreenLight, a large-scale greenhouse tomato production model that exhibits stiff dynamics.
The simulation results show that trajectory-selection RL--MPC with a one-hour prediction horizon matches the closed-loop performance of a high-performing guiding policy while significantly improving over standalone MPC with the same horizon.
On 25 weather trajectories held out during policy training, trajectory-selection RL--MPC achieved 54\% higher closed-loop performance than standalone RL and 80\% higher performance than MPC with the same horizon length.
\end{abstract}

\begin{graphicalabstract}\centering
\includegraphics[width=\textwidth]{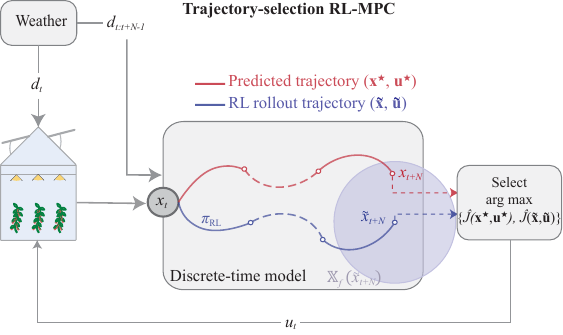}
\end{graphicalabstract}

\begin{highlights}
\item Trajectory-selection RL--MPC improves short-horizon MPC for greenhouse production control
\item Applied to GreenLight, a large-scale nonlinear greenhouse tomato production model
\item RL policy rollouts incorporate long-term economic information into short-horizon MPC
\item Improved closed-loop performance over RL policy on unseen weather trajectories
\item The source code is publicly available at \url{https://github.com/BartvLaatum/GL-Gym-MPC}
\end{highlights}


\maketitle

\section{Introduction}\label{sec:introduction}
Changing climate conditions and more extreme weather events threaten the stable production of fruits and vegetables in open-field farming \citep{silvaAssessingImpactGlobal2017, bisbisPotentialImpactsClimate2018,muldersExtremeDroughtRainfall2024}.
Greenhouse crop production systems can provide more stable food production by protecting crops against unfavorable external weather conditions \citep{faoGoodAgriculturalPratices2013}.
However, greenhouse economic performance remains strongly influenced by external weather through its effects on indoor climate, resource use, and operating costs \citep{fanourakisClimateChangeImpacts2025}.
As a result, greenhouse climate control is not only a climate management problem, but also a long-horizon economic decision problem that must balance crop yield and quality against operating cost over long horizons under exogenous weather disturbances.
In practice, greenhouse growers need to make decisions over long production horizons while accounting for economic performance indicators, future weather conditions, and crop-development dynamics \citep{vanstratenOptimalGreenhouseCultivation2010}.
The increasing scarcity of experienced employees with horticultural knowledge further motivates the need for automated control systems that can support growers in maintaining optimal growing conditions for crops at minimal cost \citep{ariesen-verschuurDigitalTwinsGreenhouse2022}.
Advanced automated control methods, such as optimal control and reinforcement learning, provide frameworks for addressing this economic control problem.

Greenhouse crop production control has been widely studied as an economic optimal control problem \citep{vanhentenGreenhouseClimateManagement1994, tapEconomicsbasedOptimalControl2000, blascoModelbasedPredictiveControl2007,vanstratenOptimalGreenhouseCultivation2010, gruberNonlinearMPCBased2011, itoGreenhouseTemperatureControl2012, ramirez-ariasMultiobjectiveHierarchicalControl2012, montoyaHybridcontrolledApproachMaintaining2016, kuijpersLightingSystemsStrategies2021, finkComparisonDynamicTomato2023}.
In this problem formulation, a controller optimizes input trajectories to maximize economic return over a given time horizon subject to system constraints under exogenous disturbances.
In greenhouse crop production systems, climate dynamics have significantly faster timescales (minutes) compared to crop development dynamics (weeks).
For fruit-producing greenhouse systems, the current operating costs for controlling the climate must be balanced against the delayed crop response times.
Ideally, the controller's optimization horizon corresponds to the length of the crop's slow timescale to balance current operating costs with the delayed benefit of fruit yield over the production period.
These delayed benefits involve fruit-development timescales of 40 to 60 days \citep{dekoningDevelopmentDryMatter1994, adamsEffectTemperatureGrowth2001}.
Solving optimal control problems online with such horizon lengths is computationally intractable \citep{finkComparisonDynamicTomato2023,xuRulebasedYearroundModel2025}.
Especially, when using greenhouse fruit-production models, such as GreenLight \citep{katzinGreenLightOpenSource2020}, which are large, nonlinear, and exhibit stiff dynamics, further challenging the numerical solution of these long-horizon optimal control problems.
This motivates the use of approximate solution methods for greenhouse climate control.

Model predictive control (MPC) is a specific optimal control method that approximates the long-horizon economic optimal control problem by optimizing over a finite prediction horizon.
At each time step, MPC optimizes the future input trajectory from the current state, applies the first optimized control input, and repeats this procedure at the next time step in a receding-horizon manner.
Since the optimization is solved online, MPC can explicitly account for predicted future disturbances such as external weather.
However, when MPC employs large-scale nonlinear prediction models rather than typical control-oriented greenhouse models, prediction horizons that capture the relevant crop-development timescale become computationally infeasible.
GreenLight \citep{katzinGreenLightOpenSource2020}, for instance, uses 27 states to model the greenhouse system, whereas the control-oriented Van Henten model \citep{vanhentenGreenhouseClimateManagement1994} is a four-state model.
As a result, practical MPC implementations for greenhouse crop production using large-scale models may rely on relatively short prediction horizons.
Such short horizons are unable to balance current operating costs with the delayed economic benefit of crop yield over the production period, resulting in myopic controllers.
These challenges motivate control approaches that enhance the economic performance of MPC controllers with prediction horizons that do not span the full production cycle.

One line of research proposes a timescale decomposition of the optimal control problem  \citep{vanhentenGreenhouseClimateManagement1994, vanstratenOptimalControlGreenhouse2010,tapEconomicsbasedOptimalControl2000,VANHENTEN2009timescale,xuAdaptiveTwoTimescale2018}.
This approach first solves the optimization problem over the entire growing season for slow crop states, resulting in an open-loop trajectory for the crop state and co-state.
Next, this crop state and co-state serve as economic performance indicators in a finite-horizon optimization problem with short prediction horizons (order of hours).
However, extending this timescale decomposition to large, nonlinear, and stiff greenhouse models that account for delayed fruit development can become computationally challenging.
More recently, \citet{kuijpersFruitDevelopmentModelling2021} addressed these differing timescales by modifying the underlying fruit development model that captures the delayed benefits from fruit harvest (order of weeks).
This modification reduced the prediction horizon of a receding horizon controller to three days.
However, this method relies on the assumption that the crop has reached a fully generative state, making it inapplicable for optimizing the initial months of the generative phase.

Reinforcement learning (RL) provides an alternative approach for approximating a controller with economic performance over longer horizons for greenhouse production systems. 
RL learns a control policy (i.e., controller) from numerous system interactions during the training phase \citep{suttonReinforcementLearningIntroduction2018a}.
In greenhouse applications, this training is performed offline using a simulation model, which shifts computational effort away from online deployment \citep{morcegoReinforcementLearningModel2023, vanlaatum2025greenlight}.
Therefore, RL can be attractive when long-term economic performance is important, but solving the full optimal control problem online is computationally intractable.
However, unlike MPC, RL does not explicitly enforce system constraints.
Instead, these constraints are typically incorporated through penalty terms in the objective function. 
This approach requires careful tuning to balance economic performance against constraint satisfaction \citep{morcegoReinforcementLearningModel2023, vanlaatum2025greenlight}.
A second limitation is that RL policies may generalize poorly to circumstances not encountered during training \citep{zhangDissectionOverfittingGeneralization2018,zhangStudyOverfittingDeep2018,kirkSurveyZeroshotGeneralisation2023}.
In greenhouse control settings, RL policies may overfit to the disturbance distributions encountered during training, for instance, a narrow set of weather trajectories, which can lead to decreased performance under unseen weather conditions.
While training RL policies on broader distributions of historical weather trajectories can improve overall performance, it does not completely prevent degraded performance under unseen weather conditions \citep{zhangStudyOverfittingDeep2018,kirkSurveyZeroshotGeneralisation2023, mediratta2024generalizationgapofflinereinforcement}.
These limitations motivate combining RL with a complementary control method, such as MPC.
In such a framework, the RL policy can provide a learned approximation of longer-term economic performance beyond the finite prediction horizon, while the MPC optimization can refine the control input trajectory over a finite horizon using actual weather predictions.

These complementary strengths of MPC and RL have motivated their integration in greenhouse crop production control studies \citep{MSAAD2025449, mallickReinforcementLearningbasedModel2025,vanlaatumStochasticModelPredictive2026}.
Previous greenhouse production control studies have combined these methods in two main ways.
One line of work used a parameterized MPC controller as an RL policy and cost function approximator \citep{mallickReinforcementLearningbasedModel2025}.
By using RL to learn the parameterization of the MPC framework, the controller improved performance relative to standalone MPC and RL without substantially increasing computation. 
Another line of work used a trained RL policy to guide the MPC optimization \citep{MSAAD2025449, vanlaatumStochasticModelPredictive2026}.
The RL policy provided terminal constraints and terminal costs to the optimization problem through a policy rollout.
This policy-guided MPC approach achieved performance comparable to MPC while reducing computation time \citep{MSAAD2025449, vanlaatumStochasticModelPredictive2026}.
Beyond greenhouse control, similar ideas have been explored in the general context of nonlinear control.
\citet{ghezzi2026rolloutoptimizeonestepnewton} initialized a real-time, iteration-based MPC with an RL rollout and refined it through a single Riccati-based Newton step. 
They proved that the suboptimality of the RL rollout is reduced quadratically by one Newton step.
The authors of \citet{reiter2025ac4mpc} introduced a parallel architecture in which the MPC trajectory and the RL rollout are evaluated at each time step using a learned critic, and the lower-cost trajectory is applied to the system.
These results demonstrate the promise of combining RL and MPC in a unified framework.
However, existing greenhouse production control studies have tested such controllers on relatively small (i.e., four-state) lettuce greenhouse models and under limited weather variation.
Therefore, it remains unclear whether policy-guided MPC approaches scale to large, nonlinear, and stiff fruit-production greenhouse models, and whether MPC optimization further improves performance relative to the guiding RL policy under weather trajectories not encountered during policy training.

To address these gaps, this work introduces \textit{trajectory-selection RL--MPC} for greenhouse fruit production control.
The proposed framework extends the policy-guided MPC with an online trajectory-selection mechanism and is evaluated on GreenLight, a large, nonlinear greenhouse tomato production model with stiff dynamics, under a training dataset and an unseen test dataset of weather trajectories.
Following earlier policy-guided MPC work, the finite-horizon optimization problem incorporates terminal region constraints and terminal costs provided by a learned RL policy via rollout trajectories \citep{MSAAD2025449}.
The terminal region constraints and terminal costs serve as the mechanism to incorporate the delayed economic benefit of crop yield into the MPC's finite-horizon optimization.
Moreover, the MPC optimization locally refines the suggested RL rollout trajectories by optimizing over a finite prediction horizon using actual weather predictions, with the aim of improving performance over the standalone RL policy in weather conditions not encountered during training.
In addition to the existing framework, trajectory-selection RL--MPC evaluates the objective for both the MPC-predicted trajectory and the rollout trajectory over the prediction horizon, selects the trajectory with the best value, and applies the first input to the system.
Unlike the parallel architecture of \citep{reiter2025ac4mpc}, which evaluates trajectories using a learned critic and a policy rollout, the proposed selection mechanism evaluates both candidate trajectories directly using the finite-horizon planning objective.
This selection mechanism is designed to avoid a decrease in performance relative to the rollout when policy-guided MPC optimization does not improve this objective.
Trajectory-selection RL--MPC is evaluated using the large-scale greenhouse tomato model GreenLight \citep{katzinGreenLightOpenSource2020}, which contains slow fruit development processes.
The GreenLight model is used for both closed-loop simulation and as the prediction model in the finite-horizon optimization problem.
Since solving the control problem over the entire simulation period is computationally intractable, this study focuses on relatively short prediction horizons (one hour) to assess whether trajectory-selection RL--MPC can improve performance when computational demands limit horizon lengths.

The closed-loop performance of trajectory-selection RL--MPC is evaluated against standalone RL and MPC.
The three methods are compared in two complementary simulation experiments:
\begin{enumerate}
    \item \textit{Trajectory-selection RL--MPC guided by a specialist RL policy.} In this setting, the RL policy is trained on a single weather trajectory, and the three control methods are evaluated on that specific weather trajectory.
    \item \textit{Trajectory-selection RL--MPC guided by a generalist RL policy.} In this setting, the RL policy is trained on a broader distribution of weather trajectories, and the three control methods are evaluated on 25 weather trajectories held out during policy training.
\end{enumerate}
The specialist-guided and generalist-guided controllers are evaluated separately because they are trained for different purposes and are intended to demonstrate different controller characteristics.
The specialist-policy experiment assesses whether trajectory-selection RL--MPC with the same prediction horizon as MPC can achieve the closed-loop performance of a high-performing RL policy by transferring the policy's longer economic foresight.
The generalist-policy experiment tests whether trajectory-selection RL--MPC improves closed-loop performance over the guiding RL policy when evaluated on weather trajectories held out during training.
Because closed-loop evaluation with GreenLight as both the simulation model and the MPC prediction model is computationally demanding, the simulation experiments evaluate controller performance over closed-loop simulation periods spanning one day.
Therefore, the simulation results assess whether a guiding RL policy can incorporate longer-term economic information into the short-horizon MPC within the trajectory-selection RL--MPC framework.
The extrapolation of these results to full production cycles is discussed in Section~\ref{sec:discussion}.

Specifically, the contributions of this work are threefold:
\begin{itemize}
    \item \textit{Trajectory-selection RL--MPC}, a new policy-guided MPC framework that compares the MPC-predicted trajectory with the RL rollout trajectory and applies the first control input from the trajectory with the best finite-horizon objective value.
    The closed-loop performance of trajectory-selection RL--MPC guided by a specialist or generalist policy was evaluated and compared with standalone RL and MPC.
    \item Simulation results showing how the framework can incorporate the delayed benefit of fruit yield into an MPC optimization problem using relatively short prediction horizons for greenhouse fruit-production systems.
    In addition, the short-horizon MPC optimization refines control input trajectories using short-term weather disturbances, improving performance relative to the guiding RL policy by reducing constraint violations and operating costs when evaluated on weather trajectories not encountered during policy training.
    \item Application of the framework to GreenLight, a large, nonlinear greenhouse tomato model that exhibits stiff dynamics.
\end{itemize}

The remainder of this paper is organized as follows: Section~\ref{sec:material-methods} formulates the economic greenhouse control problem, describes the employed GreenLight model, and presents the RL, MPC, and trajectory-selection RL--MPC control methods.
Section~\ref{sec:results} presents the simulation results, and Section~\ref{sec:discussion} discusses these results.
Finally, Section~\ref{sec:conclusion} presents the paper's conclusions and outlines directions for future work.

\section{Materials and Methods}\label{sec:material-methods}
This work compared the closed-loop performance of the trajectory-selection RL--MPC framework for greenhouse tomato production control against standalone MPC and RL.
Section~\ref{sec:optimal-control-problem} formulates the economic greenhouse control problem studied in this work, after which the reward function and prediction model are specified in Sections~\ref{sec:reward} and~\ref{sec:model}.
Section~\ref{sec:rlmpc} introduces the trajectory-selection RL--MPC framework, which combines the RL policy with the MPC optimization.
The performance metrics used to evaluate these three control methods are described in Section~\ref{sec:performance-metrics}, and the weather datasets used for training and evaluation are presented in Section~\ref{sec:weather-data}.
Finally, Section~\ref{sec:simulation-experiments} describes the simulation experiments used to compare closed-loop controller performance.

\subsection{Economic greenhouse control problem}\label{sec:optimal-control-problem}
This work focuses on the economic control of a greenhouse tomato production system by maximizing its closed-loop economic performance under varying weather trajectories.
The greenhouse system is affected by controllable input $u\in\mathbb{R}^{n_{u}}$, and uncontrollable input $d\in\mathrm{\mathbb{R}^{n_{d}}}$.
The weather trajectory is denoted by $\mathbf{d}\coloneq \{d(0),\dots,d(T-1)\}$. 
For a given $\mathbf{d}$, the following scenario-conditioned economic control problem is defined as:
\begin{subequations}\label{eq:optimization_problem}
	\begin{align}
        \underset{\substack{\bigl(x(0),\, \dots,\, x(T)\bigr)\\ \bigl(u(0),\, \dots,\, u(T-1)\bigl)}}{\max} \quad & 
		\sum_{t=0}^{T-1}r\left(x(t), x(t+1), u(t)\right) \label{eq:optimization_problem:objective_function} \\
        \text{s.t.} \quad
        & \left.\begin{aligned}		
    		& x(t+1) = f\bigl(x(t),u(t),d(t)\bigr) \\
        	& u_\text{min} \leq u(t) \leq u_\text{max} \\
    		& \mspace{-4mu} \left| u(t) - u(t-1) \right| \leq \delta u_\text{max} \\
            \end{aligned}\quad\right\} \quad\forall t\in\{0,\dots, T-1\}
    \label{eq:optimization_problem:constraints}
	\end{align}
\end{subequations}
where $t\in\mathbb{Z}^{0+}$ represents the discrete time step and $T$ the end of the closed-loop period.
Let $x(t)\in\mathbb{R}^{n_{x}}$ denote the state at time step $t$, and let $r(\cdot)$ denote the economically based reward, including a penalty for violations of climatic state constraints.
The reward function is further detailed in Section~\ref{sec:reward}.
The nonlinear function $f(\cdot)$ represents the dynamical model that describes the system's next state given controllable and uncontrollable inputs. 
At each time step $t$, the employed controller computes a control input $u(t)$, which is applied to the system; iterating this procedure over $t=0, \dots,T$ yields the closed-loop state $\bigl(x(0),\, \dots,\, x(T)\bigr)$ and input trajectory $\bigl(u(0),\, \dots,\, u(T-1)\bigl)$.
In addition, the control input was subject to amplitude bounds and rate-of-change constraints.
Following previous work, to limit short-term variability in the indoor climate while ensuring smooth actuation, the control input $u(t)$ was constrained to limit the rate of change within $10\%$ ($\delta = 0.1$) of the corresponding upper boundary $u_{\max}$ \citep{boersmaRobustSamplebasedModel2022,morcegoReinforcementLearningModel2023,svensenChanceconstrainedStochasticMPC2024,mallickReinforcementLearningbasedModel2025,MSAAD2025449,vanlaatum2025greenlight,vanlaatumStochasticModelPredictive2026}.

In this work, $\mathbf{d}$ was selected from a finite dataset of weather trajectories, i.e., $\mathbf{d}\in\mathcal{D}^{(i)}$, where $\mathcal{D}^{(i)}$ denotes the dataset used in a given simulation experiment, e.g., for training an RL policy or evaluating controller performance ($i\in\{\mathrm{train,\, test}\}$).
The weather sets used in the simulation experiments are defined in Section~\ref{sec:weather-data}.

\subsection{Tomato greenhouse climate model}\label{sec:model}
The nonlinear dynamical tomato greenhouse climate model considered in this work is based on GreenLight \citep{katzinGreenLightOpenSource2020}.
GreenLight is a high-fidelity tomato greenhouse climate model for simulating indoor climate in high-tech greenhouse facilities and has been validated against real-world greenhouse data.
The tomato model is a simplified version of the tomato crop-yield model of \cite{vanthoorMethodologyModelbasedGreenhouse2011}.
This simplification may influence the timing of fruit harvest, but the total yield is not expected to be affected \citep{katzinGreenLightOpenSource2020}.

The resulting state-space model comprises indoor-climate states (temperature, \ce{CO2} concentration, and humidity) and crop states representing dry weight and crop development. 
Table~\ref{tab:gl_states} lists all 27 state variables.
Throughout this work, the nonlinear system dynamics are defined by the following continuous-time greenhouse model:
\begin{equation}\label{eq:continuous-model}
   \dot{x} = f_c\bigl(x(t_c),u(t_c),d(t_c)\bigr),
\end{equation}
where $f_c(\cdot)$ describes the system of differential equations, and $t_c\in\mathbb{R}$ represents the continuous time.
To perform closed-loop simulations and apply the control methods considered in this work, a discrete-time form of this model is required.
This discrete-time model was used both to generate the closed-loop simulation experiments and within the control methods: RL interacted with it during training, while MPC-based controllers used it in the finite-horizon optimization problem.
The resulting discrete-time greenhouse model is given by:
\begin{equation}\label{eq:greenlight-model}
        \begin{aligned}
        x(t+1) &= f(x(t), u(t), d(t)),\\
        y(t)   &= h(x(t)),
        \end{aligned}
\end{equation}
where $x(t)\in\mathbb{R}^{27}$ denotes the state, $y(t)\in\mathbb{R}^{27}$ the output, $u(t)\in\mathbb{R}^{6}$ the control input, and $d(t)\in\mathbb{R}^{7}$ the weather disturbances.
The nonlinear function $f(\cdot)$ denotes the numerical integration function, and $h(\cdot)$ is the output function.
This work used a discretization time (and control interval) of $\Delta t=300\,\mathrm{s}$, i.e., five minutes.
This work considered six input variables that manipulate the indoor climate; their lower and upper bounds are listed in Table~\ref{tab:gh_input}.
For a complete model description, the reader is referred to \citet{katzinGreenLightOpenSource2020}.

\begin{table}[t]
\centering
\caption{Control input $u$ to the tomato greenhouse model~\eqref{eq:greenlight-model} and corresponding lower and upper bounds.}
\label{tab:gh_input}
\begin{tabular}{@{}lcccp{8.2cm}@{}}
\toprule
\textbf{Symbol} & \textbf{Lower bound} ($u_{\min}$) & \textbf{Upper bound} ($u_{\max}$) & \textbf{Unit} & \textbf{Description} \\
\midrule
$u_{
\mathrm{boil}}$ & $0$ & $44$ & $\mathrm{W\,m^{-2}}$ & The energy from the boiler to the heating pipes \\
$u_{\mathrm{led}}$ & $0$ & $116$ & $\mathrm{W\,m^{-2}}$ & Electrical power to LED lighting \\
$u_{\ce{CO2}}$ & $0$ & $5$ & $\mathrm{mg\,m^{-2}\,s^{-1}}$ & Greenhouse CO$_2$ injection \\
$u_{\mathrm{vent}}$ & $0$ & $1$ & -- & Opening of the ventilation vents \\
$u_{\mathrm{thscr}}$ & $0$ & $1$ & -- & Thermal screen set (1 represents fully deployed) \\
$u_{\mathrm{blscr}}$ & $0$ & $1$ & -- & Blackout screen set (1 represents fully deployed) \\
\bottomrule
\end{tabular}
\end{table}

The continuous model~\eqref{eq:continuous-model} exhibits stiff dynamics, which has implications for numerical integration of the model and controller design.
Obtaining stable solutions of the discrete-time form of the model~\eqref{eq:greenlight-model} requires small integration step sizes ($1\mathrm{s}$) when using explicit solvers, such as n\textsuperscript{th}-order Runge-Kutta methods.
Therefore, this work used CVODES \citep{serban2005cvodes}, an implicit, variable-order, variable step-size integration method with automatic differentiation, to obtain a numerical solution to~\eqref{eq:greenlight-model}.
Although accurate, this integration method is computationally demanding when solving nonlinear finite-horizon optimization problems. 
Nonetheless, this approach proved more efficient than explicit methods with sufficiently small step sizes.
Because solving the economic control problem in~\eqref{eq:optimization_problem} over the full horizon $T$ is computationally intractable, approximate solution methods were considered instead.
Specifically, this work studies RL, MPC, and a combined RL--MPC approach, which are introduced next.

\subsection{Reward function}\label{sec:reward}
The reward function $r(\cdot)$ defining the economic control problem in~\eqref{eq:optimization_problem} comprises two terms: an economic-based term and a penalty term that quantifies constraint violations.
Combining the economic reward $r_{\mathrm{e}}(\cdot)$ with the penalty function $r_{\mathrm{p}}(\cdot)$ yields the following reward function:
\begin{equation}\label{eq:reward}
    r\left(x(t), x(t+1), u(t)\right) =r_{\mathrm{e}}\bigl(x(t), x(t+1), u(t) \bigr) - r_{\mathrm{p}}\bigl(x(t+1)\bigr),
\end{equation}
where this definition is chosen such that maximizing $r(\cdot)$ corresponds to maximizing economic performance while minimizing constraint violations.
e
The economic reward is defined by two terms: revenue from tomato yield and operating costs associated with the controlled input, i.e., heating, \ce{CO2} injection, and lighting electricity.
To quantify the economic performance and enable optimization by the employed controller, the following economic reward function is defined:
\begin{equation}\label{eq:economic-reward}
    \begin{aligned}
        r_{\mathrm{e}}\left(x(t), x(t+1), u(t)\right) =  & ~c_{\mathrm{frt}}\left(x_{\mathrm{frt}}(t+1)-x_\mathrm{frt}(t)\right)-\\ 
        &\left(c_{\mathrm{CO_2}}u_{\mathrm{CO_2}}(t) + c_{\mathrm{heat}}u_{\mathrm{heat}}(t) + c_{\mathrm{elec}}u_{\mathrm{led}}(t)\right)\Delta t,
    \end{aligned}
\end{equation}
where $c_{\mathrm{frt}}$ is the unit selling price of the tomato fruits, $x_{\mathrm{frt}}(t)$ is the amount of fruit yield at time $t$, $c_{i}$ represent the unit costs of of input $u_{i}$, and $\Delta t$ is the control interval.
The values of the unit price and costs $c_{i}$ for the simulation experiments are defined in Section~\ref{sec:simulation-experiments}.

Practical greenhouse operation also requires maintaining indoor climate conditions within acceptable ranges.
Unfavorable indoor climate conditions can induce crop diseases, which negatively affect fruit yield.
Because the model does not capture this relationship between crop disease and yield, it does not directly affect the economic reward function.
To better reflect greenhouse practice, the optimization problem includes constraints on three indoor climate variables: temperature, humidity, and \ce{CO2} concentration denoted by $(y_{\mathrm{T,air}}, y_{\mathrm{CO_2,air}}$, and $y_{\mathrm{RH,air}})$, respectively.
As in previous work, these constraints were implemented as soft constraints by penalizing violations of the constraints with a linear penalty function \citep{morcegoReinforcementLearningModel2023,MSAAD2025449,mallickReinforcementLearningbasedModel2025,vanlaatumStochasticModelPredictive2026}.
This penalty function is defined as:
\begin{equation}\label{eq:penalty-cost}
    r_{\mathrm{p}}(x(t)) = g_{\Tair}\bigl(y_{\Tair}(t)\bigr) + g_{\COair}\bigl(y_{\COair}(t)\bigr)+g_{\mathrm{RH, air}}\bigl(y_{\mathrm{RH, air}}(t)\bigr),
\end{equation}
where the outputs $y_{i}(t)$ are obtained from the output function in~\eqref{eq:greenlight-model}.
The functions $g_{i}$ apply penalties for deviations from their \emph{a priori} defined operating ranges.
Each penalty function is defined as:
\begin{equation}\label{eq:penalty_factors}
	g_{i}\bigl(y_{i}(t)\bigr) = 
	\begin{cases}
		\lambda_{i} \bigl(y_{i}(t) - y_{i}^\text{max}\bigr) & \text{if } y_{i}(t) > y_{i}^\text{max}, \\
		\lambda_{i} \bigl(y_{i}^\text{min} - y_{i}(t)\bigr) & \text{if } y_{i}(t) < y_{i}^\text{min}, \\
		0 & \text{otherwise},
	\end{cases}
\end{equation}
where $\lambda_{i}$ denotes the penalty coefficients that regulate the scale of each penalty term.
The choice of these coefficients strongly affects the controller's performance, as larger values emphasize constraint violations over economic performance, whereas smaller values do the opposite.
The penalty coefficients were chosen such that the penalty term was sufficiently large to discourage constraint violations, yet remained of the same order of magnitude as the economic reward in~\eqref{eq:economic-reward}.
This scaling prevented the penalty term from dominating the objective while providing a practical balance between economic performance and constraint satisfaction.
Section~\ref{sec:simulation-experiments} defines the values of these penalty coefficients $\lambda_{i}$.

\subsection{RL}\label{sec:rl}
This work leverages the Proximal Policy Optimization (PPO) algorithm \citep{schulmanProximalPolicyOptimization2017} to train an RL policy, denoted by $\piRL$.
Two PPO-based policies are considered: a \emph{specialist policy} and a \emph{generalist policy}.
Both policies share the same training objective and action-to-input mapping, but differ in their policy input, hyperparameter settings, and training weather data.

Each policy is trained by interacting with the dynamical model in~\eqref{eq:greenlight-model} using historical weather trajectories $\mathbf{d}$ sampled from the training dataset.
The goal is to approximate the economic control problem in~\eqref{eq:optimization_problem} by optimizing the expected discounted return:
\begin{equation}\label{eq:rl-objective}
    \max_{\pi_{\mathrm{RL}}}\mathbb{E}_{(\mathbf{x},\mathbf{u})\sim\rho(\pi_{RL}, \mathbf{d})}
    \left[
    \sum_{k=t}^{T-1} \gamma^{k-t}r\bigl(x(k), x(k+1), u(k)\bigr)
    \right],
\end{equation}
where $x(k)$ is subject to the dynamics in~\eqref{eq:greenlight-model}.
and control input is sampled from $u(k)\sim\piRL(\cdot)$.
The joint distribution $\rho$ is induced by the policy and the sampled weather trajectory $\mathbf{d}$.
At $k=T$, the system transitions to an absorbing terminal state that yields zero reward indefinitely, allowing us to write the infinite sum of rewards using a finite-time horizon \citep{suttonReinforcementLearningIntroduction2018a}.
Ideally, RL would directly optimize the undiscounted objective posed in~\eqref{eq:optimization_problem}.
However, computing undiscounted returns over long horizons introduces high variance in gradient estimates, which can cause convergence issues when training policies using deep neural networks.
To address this, a discount factor of $\gamma=0.995$ is applied to limit this variance.
Given a five-minute control interval ($\Delta t=300\mathrm{(s)}$), this corresponds to an effective horizon of approximately $\frac{1}{1-\gamma}=200$ time steps, i.e., roughly 16 hours \citep{kearns2002near}.
Although this discounted objective acts as a surrogate for the objective posed in~\eqref{eq:optimization_problem}, it represents a necessary bias-variance trade-off to ensure the stable convergence of a policy that performs well on the undiscounted closed-loop objective.


At each time step $t$, the policy receives a feature vector $\phi(t)$ and produces an action $a(t)=\piRL(\phi(t))$. 
This action $a(t)$ is then mapped to the controllable input $u(t)$:
\begin{equation}
    u(t) = \max\Big(~u_{\min},~\min\big(~u(t-1) + a(t)\delta u_{\max}~, u_{\max}\big)~\Big).
\end{equation}
Hence, the action scales the input rate of change relative to the previously applied control input. 
Because $a(t)\in[-1, 1]$, the input rate constraints are satisfied by construction. 
To comply with the input bounds, the resulting control input is clipped between $u_{\min}$ and $u_{\max}$.

The policy feature vector $\phi(t)$ is constructed from a common base observation vector $o(t)$.
At each time step $t$, this vector is defined as:
\begin{equation}
\label{eq:base_observation}
\begin{aligned}
o(t)
\;&=\;
\Big[
\begin{matrix}
x_{\mathrm{T,air}}(t) &
x_{\ce{CO2}}(t) &
        x_{\mathrm{RH}}(t)&
        x_{\mathrm{T,pipe}}(t)&
        x_{\mathrm{frt}}(t)& \quad \dots
\end{matrix}
\\
&\qquad \qquad \dots
\quad \begin{matrix}
\sin\!\big(2\pi\,\tau_{\mathrm{day}}(t)\big)&
\cos\!\big(2\pi\,\tau_{\mathrm{day}}(t)\big)&
\sin\!\big(2\pi\,\tau_{\mathrm{year}}(t)\big)&
\cos\!\big(2\pi\,\tau_{\mathrm{year}}(t)\big)& \quad \dots
\end{matrix}
\\
&\qquad \qquad \qquad \qquad \qquad \qquad \dots
\quad \begin{matrix}
u(t-1)&
d_{\mathrm{iGlob}}(t)&
d_{\mathrm{T}}(t)&
d_{\ce{CO2}}(t)&
d_{\mathrm{RH}}(t)&
d_{\mathrm{wind}}(t)
\end{matrix}
\Big] \in \mathbb{R}^{\,20} 
\end{aligned}
\end{equation}
where $u(t-1)\in\mathbb{R}^{6}$ denotes the previously applied control input,
$\tau_{\mathrm{day}}(t)\in[0,1)$ the normalized time of day, and
$\tau_{\mathrm{year}}(t)\in[0,1)$ the normalized day of year.
Outdoor soil and sky temperature are excluded, as these are not commonly measured in commercial greenhouse facilities \citep{petropoulou2023lettuce, maree2025autonomous}.

For the specialist policy, short-term temporal context is incorporated by stacking the current and previous base observations:
\begin{equation}\label{eq:specialist-input}
    \phi_{\mathrm{specialist}}(t)
    \;=\;
    \big[o(t-1)^\top\; o(t)^\top\big]^\top.
\end{equation}
For the generalist policy, the observation is augmented with the weather disturbances over the next 12 time steps (one hour):
\begin{equation}
\label{eq:observation_vector_forecast}
\begin{aligned}
o_{\mathrm{fc}}(t)
\;&=\;
\Big[
\begin{matrix}
o(t)^\top &
d_{\mathrm{iGlob}}(t{+}1{:}t{+}12)^\top &
d_{\mathrm{T}}(t{+}1{:}t{+}12)^\top \quad \dots
\end{matrix}
\\
&\qquad \qquad \qquad \qquad \dots
\begin{matrix}
\quad d_{\ce{CO2}}(t{+}1{:}t{+}12)^\top &
d_{\mathrm{RH}}(t{+}1{:}t{+}12)^\top &
d_{\mathrm{wind}}(t{+}1{:}t{+}12)^\top
\end{matrix}
\Big] \in \mathbb{R}^{\,80} 
\end{aligned}
\end{equation}
where \[d_{i}(t{+}1{:}t{+}12)\triangleq\big[d_{i}(t{+}1),\ldots,d_{i}(t{+}12)\big]^\top\] 
stacks the next 12 values of the weather variable $i$.
Similar to the specialist policy, the historic temporal context is incorporated by stacking the current and the previous observations:
\begin{equation}\label{eq:generalist-input}
\phi_{\mathrm{generalist}}(t)
    \;=\;
    \big[o_{\mathrm{fc}}(t-1)^\top\;o_{\mathrm{fc}}(t)^\top\big]^\top.
\end{equation}
To isolate the effects of the trajectory-selection RL--MPC framework structure and weather-disturbance variability on controller performance, both policies are provided with uncertainty-free feature vectors $\phi(t)$.
In particular, no uncertainty in the weather forecast is considered in $\phi_{\mathrm{generalist}}(t)$.

A random search informed the PPO hyperparameters for both policies. 
Because the generalist policy has a higher-dimensional input space, its hyperparameters require additional tuning.
The resulting hyperparameter settings are listed in Table~\ref{tab:ppo_hyperparams}.

The specialist policy was developed to obtain a high-performing policy for a single scenario by training exclusively on a single weather trajectory, to assess whether trajectory-selection RL--MPC can maintain the specialist policy's performance.
The generalist policy was used to study how a policy trained across multiple weather trajectories can improve the performance of trajectory-selection RL--MPC under unseen weather conditions.
Section~\ref{sec:weather-data} describes the training dataset used for each policy in more detail.
The RL environment, including observations, rewards, and weather-trajectory sampling, was implemented in GreenLight-Gym \citep{vanlaatum2025greenlight}.
Policy training was performed using Stable-Baselines-3 \citep{raffinStableBaselines3ReliableReinforcement2021}.




\begin{table}[t]
\centering
\footnotesize
\setlength{\tabcolsep}{5pt}
\caption{PPO hyperparameters used for training the specialist and generalist policies. 
The generalist policy differs from the specialist policy only in the batch size, learning rate, and number of environment steps per update ($n_{\text{steps}}$).}
\label{tab:ppo_hyperparams}
\begin{tabular}{@{}lcc@{}}
\toprule
\textbf{Hyperparameter} & \textbf{Specialist} & \textbf{Generalist} \\
\midrule
Total time steps & $4{,}000{,}000$ & $4{,}000{,}000$ \\
Number of parallel envs ($n_{\text{envs}}$) & $8$ & $8$ \\
Epochs per update ($n_{\text{epochs}}$) & $8$ & $8$ \\
Discount factor ($\gamma$) & $0.995$ & $0.995$ \\
GAE parameter ($\lambda$) & $0.95$ & $0.95$ \\
Clipping range & $0.2$ & $0.2$ \\
Entropy coefficient ($c_{\text{ent}}$) & $1.8\times 10^{-4}$ & $1.8\times 10^{-4}$ \\
Value function coefficient ($c_{\text{vf}}$) & $0.896$ & $0.896$ \\
Policy network size ($\pi$) & $[256,\,256]$ & $[256,\,256]$ \\
Value network size ($V_{\theta}$) & $[256,\,256]$ & $[256,\,256]$ \\
Activation function & \texttt{SiLU} & \texttt{SiLU} \\
\midrule
\textit{Learning rate} & $8.75\times 10^{-5}$ & $5.0\times 10^{-5}$ \\
\textit{environment steps per update} ($n_{\text{steps}}$) & $1024$ & $512$ \\
\textit{Batch size} & $128$ & $64$ \\
\bottomrule
\end{tabular}
\end{table}

\subsection{MPC}\label{sec:mpc}
Given the reward function in~\eqref{eq:reward}, and the discrete-time model in~\eqref{eq:greenlight-model}, at each time step $t$ the MPC is provided with the current state $x(t)$ and the previously applied input $u(t-1)$, and solves the following finite horizon optimization problem:
\begin{subequations}\label{eq:MPC_optimization_problem}
	\begin{align}
		\underset{\mathbf{x}, \mathbf{u}}{\max}
		  \quad &\sum_{k=t}^{t+N-1}r\bigl(x(k), x(k+1), u(k)\bigr) \\
		\text{s.t.} \quad
        & \left.
        \begin{aligned}            
    		& x(k+1) = f(x(k),u(k),d(k)) \\
    		& u_\text{min} \leq u(k) \leq u_\text{max} \\
    		& \left|u(k) - u(k-1)\right| \leq \delta u_\text{max}
        \end{aligned}
        \quad\right\}\;\; \forall k \in \{t, \dots, t+N-1\}
\end{align}
\end{subequations}
where $x(t)$ and $u(t-1)$ enter the problem as fixed initial conditions.

The solution of the optimization problem is denoted by the state and input trajectories:
\begin{equation}
\begin{aligned}    
    \mathbf{x}^\star &= \bigl(x^\star(t),\, \dots, \,x^\star(t+N)\bigr), \quad
    \mathbf{u}^\star &= \bigl(u^\star(t),\, \dots, \,u^\star(t+N-1)\bigr).
\end{aligned}
\end{equation}
From the computed input trajectory $\mathbf{u}^\star$, only the first element is applied to the system; that is, $u(t)=u^\star(t)$.
The optimization problem is solved in a receding-horizon fashion.
At time step $t+1$, the problem is solved again using the updated initial conditions $x(t+1)$ and $u(t)$, and the procedure is repeated.
At each time step, the MPC is assumed to have access to full state information and perfect weather forecasts.

Due to the variable-step integrator required to evaluate~\eqref{eq:greenlight-model}, the computational demand of solving this optimization problem is high, and the MPC is therefore restricted to relatively short prediction horizons (e.g., 1 hour) that do not capture the entire closed-loop simulation period ($N << T$).
The prediction horizon used in the simulation experiments are specified in Section~\ref{sec:simulation-experiments}.

The nonlinear MPC optimization problem is implemented using CasADi \citep{anderssonCasADiSoftwareFramework2019a} and solved with the interior-point optimization software IPOPT \citep{wachterImplementationInteriorpointFilter2006}.
Several IPOPT settings were adjusted to make the optimization computationally tractable.
Specifically, the Hessian is computed using a quasi-Newton approximation.
In addition, the feasibility and acceptable tolerances are relaxed by setting the constraint violation tolerance to 0.01, the acceptable constraint violation tolerance to 0.1, and the acceptable optimality tolerance to 0.1.
These settings allow IPOPT to terminate at approximately feasible and approximately (locally) optimal solutions, thereby reducing computational demand and improving convergence.\footnote{Note that local optimality is guaranteed only up to some tolerance, which is a typical (but not always discussed) aspect of numerical optimization for MPC. Hence, other trajectories may exist that achieve better performance than this computed solution, but locally optimal solutions are often sufficient for good MPC performance in industrial practice.} 

\subsection{Trajectory-selection RL--MPC}\label{sec:rlmpc}
The proposed trajectory-selection RL--MPC framework extends policy-guided MPC, in which a trained RL policy provides terminal ingredients for a finite-horizon MPC optimization problem \citep{reiter2025ac4mpc, MSAAD2025449}.
The trained RL policy is used to provide longer-term economic information into short-horizon MPC by defining a terminal region constraint and terminal cost.
The MPC optimization can refine the control input trajectory by solving a finite-horizon optimization problem using the prediction model and actual weather disturbances.
The RL rollout trajectory and the MPC-predicted trajectory are evaluated on the finite-horizon planning objective, and the trajectory with the better objective value is selected.
This trajectory selection mechanism is important, as the computed solution from the finite-horizon optimization problem is not guaranteed to outperform the RL rollout. 
This observed behavior may occur due to the presence of local minima and the relatively high tolerance required for computationally tractable solutions (as stated in the previous section) to this optimization problem. 
Through this mechanism, we aim to avoid the possibility that the proposed framework degrades the performance of the original RL policy due to numerical or computational issues.
Together, these components are designed to improve the closed-loop performance of short-horizon MPC controllers for greenhouse tomato production across weather conditions not encountered during RL policy training.

\begin{figure}
    \centering
    \includegraphics[width=0.5\linewidth]{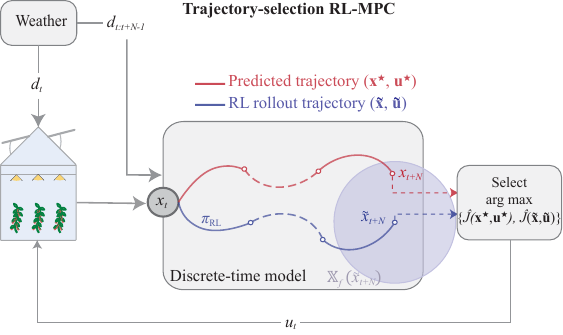}
    \caption{\textbf{Illustrative sketch of trajectory-selection RL--MPC in closed-loop simulation.}
    At time step $t$, from the measured greenhouse state $x(t)$, the RL policy generates a rollout trajectory ($\mathbf{\tilde x}, \mathbf{\tilde u}$) over the prediction horizon $N$ using future weather disturbances $d_{t:t+N-1}$.
    Next, MPC solves a finite-horizon optimization problem over the same horizon $N$, yielding the predicted trajectory ($\mathbf{x}^\star,\mathbf{u}^\star$).
    The shaded region represents the terminal region constraint $\mathbb{X}_f$, constructed around the terminal state of the RL rollout, $\tilde x(t+N)$.
    The terminal state $x^\star(t+N)$ of the MPC-predicted trajectory is constrained to end in this terminal region.
    The first control input $u(t)$ of the trajectory with the better finite-horizon objective value $\hat J$ is applied to the greenhouse system.
    Both the RL rollout and the finite-horizon optimization use the discrete-time model $f(\cdot)$~\eqref{eq:greenlight-model} to predict the evolution of the state trajectory.
    }
    \label{fig:sketch-selection-rlmpc}
\end{figure}

Figure~\ref{fig:sketch-selection-rlmpc} illustrates the closed-loop implementation of trajectory-selection RL--MPC. 
At each time step $t$, the current greenhouse state $x(t)$ is measured and a weather forecast $d(t:t+N-1)$ is provided over the prediction horizon $N$.
Starting from this state, the trained RL policy $\piRL$ generates a rollout trajectory over the same prediction horizon.
The terminal state of this rollout is then used to construct a terminal region constraint.
The framework then solves a finite-horizon optimization problem over the same horizon $N$, yielding a predicted trajectory, whose terminal state $x^\star(t+N)$ is constrained to lie within this terminal region.
After optimization, both the RL rollout trajectory and the MPC-predicted trajectory are evaluated using the same finite-horizon planning objective.
The trajectory with the better objective value is selected, and only the first input $u(t)$ of the selected trajectory is applied to the greenhouse system.
Similar to nominal MPC, the framework is assumed to have access to full state information and perfect weather forecasts.

The trained RL policy $\piRL$ is rolled out in simulation from the current state $x(t)$, yielding a state and input trajectory over the prediction horizon $N$:
\begin{equation}\label{eq:rl-rollout}
\mathbf {\tilde{x}}=(\tilde x(t),\, \dots, \,\tilde x(t+N)), \quad \mathbf{\tilde{u}} = (\tilde u(t),\, \dots\, \tilde u(t+N-1)).
\end{equation}
This rollout uses the same policy input $\phi(t)$ construction and action-to-input mapping as described in Section~\ref{sec:rl}.

Trajectory-selection RL--MPC incorporates a terminal region constraint and a terminal cost into the optimization problem.
The terminal region constraint is defined around the terminal state of the RL rollout $\tilde x(t+N)$ and constrains the predicted terminal state to remain within a region around this suggested terminal state:
\begin{equation}\label{eq:terminal-constraint}
    \mathbb{X}_{f}\bigl(\tilde{x}(t+N)\bigr)=\left\{ x\in\mathbb{X} \,|\, |x-\tilde{x}(t+N)| \leq \omega \bar x_{\mathrm{ref}}\right\},
 \end{equation}
where the inequality is applied element-wise, $\bar x_{\mathrm{ref}}$ is a normalization vector,  $\omega>0$ scales the region size proportional to this vector, and $\tilde x(t+N)\in \mathbb{X}$ denotes the terminal state of the RL rollout.

The terminal cost penalizes deviations from the rollout terminal state, allowing the optimizer to trade off deviations from the terminal state against improvements in the objective.
The MPC optimization can improve over the rollout trajectory while remaining close to the terminal state.
However, if the RL policy provides a poor terminal state, these terminal ingredients may restrict the improvement MPC can achieve.
The terminal cost is defined as follows:
\begin{equation}\label{eq:terminal-cost}
    \hat V_{f}\bigl(x,\tilde x(t+N)\bigr)=\rho\sum_{i=1}^{n}\left| \frac{x_i-\tilde{x}_i(t+N)}{\bar x_{\mathrm{ref},\,i}}\right|,
\end{equation}
where $\rho$ is a scaling factor and $\bar x_{\mathrm{ref}}$ normalizes the state differences element-wise.
A linear $L_{1}$ penalty is used instead of a quadratic $L_{2}$ penalty to prevent large deviations in individual states from dominating the terminal cost.
For all simulation experiments, the design parameters are set to $\omega=0.05$ and $\rho=10^{-5}$.

The normalization vector $\bar x_{\mathrm{ref}}$ is computed once at the start of each simulation.
First, a reference trajectory \[x_{\mathrm{ref}}=\big(x(0), \dots, x(T)\big)\] is generated by rolling out $\piRL$ from the initial state $x(0)$ under the weather trajectory $\mathbf{d}$ using the dynamical model in~\eqref{eq:greenlight-model}.
Next, the average state of this reference trajectory is computed as:
\begin{equation}
\bar x_{\mathrm{ref}}=\frac{1}{T+1}\sum_{k=0}^{T} x_{\mathrm{ref}}(k).
\end{equation}
This normalization affects only the scaling of the terminal region constraint and terminal cost, not the framework structure.

\subsubsection*{Finite-horizon optimization problem}
Combining the reward with the terminal cost yields the finite-horizon planning objective over the horizon $N$:
\begin{equation}\label{eq:objective-planning}
    \hat J(\mathbf x, \mathbf u,\tilde{x}(t+N)) = \sum_{k=t}^{t+N-1}r\bigl(x(k), x(k+1), u(k)\bigr) - \hat{V}_{f}\bigl(x(t+N), \tilde x(t+N)\bigr).
\end{equation}
Thus, the framework aims to maximize cumulative economic reward over the horizon while penalizing deviations from the RL rollout terminal state.

Given the planning objective $\hat J$ in~\eqref{eq:objective-planning}, the dynamical model in~\eqref{eq:greenlight-model}, and the terminal region constraint in~\eqref{eq:terminal-constraint}, trajectory-selection RL--MPC solves the following finite-horizon optimization problem:
\begin{subequations}\label{eq:RL--MPC}
\begin{align}
    \max_{\mathbf{x},\mathbf{u}} \quad
    & \hat J(\mathbf{x},\mathbf{u},\tilde{x}(t+N)) \\
    \text{s.t.} \quad
    &  \left.%
       \begin{aligned}
            & x(k+1) = f\bigl(x(k),u(k),d(k)\bigr) \\
            & u_{\min} \leq u(k) \leq u_{\max} \\
            & |u(k)-u(k-1)| \leq \delta u_{\max}
        \end{aligned}
        \qquad\right\}
    \qquad \forall k \in \{t,\dots,t+N-1\} \\
    & x(t+N) \in \mathbb{X}_f\bigl(\tilde x(t+N)\bigr)
\end{align}
\end{subequations}
where the current state $x(t)$ and the previously applied input $u(t-1)$ are set as fixed initial conditions.
The RL rollout trajectories in~\eqref{eq:rl-rollout} are used as a warm start for solving~\eqref{eq:RL--MPC}.
The trajectory-selection RL--MPC framework is assumed to have full access to the system state.
The solution of the optimization problem is denoted by the state and input trajectories:
\begin{equation}
\begin{aligned}    
    \mathbf{x}^\star &= \bigl(x^\star(t),\, \dots, \,x^\star(t+N)\bigr), \quad
    \mathbf{u}^\star &= \bigl(u^\star(t),\, \dots, \,u^\star(t+N-1)\bigr).
\end{aligned}
\end{equation}

After optimization, both the predicted trajectories ($\mathbf{x}^\star, \mathbf{u}^\star$) and the rollout trajectories ($\mathbf{\tilde x}, \mathbf{\tilde u}$) are evaluated using the finite-horizon planning objective~\eqref{eq:objective-planning}.
For the RL rollout, the terminal penalty is zero by construction, so deviations from the rollout terminal state are penalized only for the MPC-predicted solution.
Finally, the first input of the trajectory with the higher objective value is applied to the system.
This \emph{selection mechanism} allows the framework to fall back to the RL rollout when the optimization does not improve upon it.
Note that the planning objective is undiscounted, whereas the policy was trained using a discounted objective.
This choice is justified since this work aims to optimize the undiscounted closed-loop performance; thus, selecting on a discounted objective would not align with the undiscounted evaluation criterion.
As a result, the selection mechanism may favor MPC-predicted trajectories since the RL policy was optimized for a different (discounted) objective.

At time step $t+1$, the optimization problem is solved again using the updated state measurement, the new RL rollout, and the resulting new terminal region, and the procedure is repeated.
Algorithm~\ref{alg:rlmpc_closedloop} provides a complete closed-loop description.

The finite-horizon optimization problem was implemented and solved using the same software and optimization settings as the nominal MPC framework, see Section~\ref{sec:mpc}.

\begin{algorithm}[t]
\caption{Trajectory-selection RL--MPC controller in closed-loop simulation}
\label{alg:rlmpc_closedloop}
\begin{algorithmic}[1]
\Require \begin{tabular}[t]{@{}l@{\,: }l@{}}
$f(\cdot)$            & dynamical model~\eqref{eq:greenlight-model}\\
$x(0),\,u(0)$         & initial state/input\\
$\pi_{\mathrm{RL}}$   & learned policy\\
$\mathbf d$           & weather trajectory\\
$p$                   & model parameters\\
$N,\,T$               & prediction horizon, length of simulation
\end{tabular}
\State Compute $\bar x_{\mathrm{ref}}$ using $\piRL$, $\mathbf{d}$, starting at $x(0)$ and subject to $f(\cdot)$

\For{$t=0,\dots,T-1$}
    \State Obtain RL rollout over $N$: $\bigl(\mathbf{\tilde{x}},\,\mathbf{\tilde{u}}\bigr)$
    \State Compute objective of the rollout: $\hat J_{\mathrm{RL}}\gets \hat J\bigl(\mathbf{\tilde{x}},\, \mathbf{\tilde{u}}\bigr)$
    \Comment{see~\eqref{eq:objective-planning}}
    \State Construct terminal region constraint: $\mathbb X_f(\tilde x(t+N))$
    \Comment{see~\eqref{eq:terminal-constraint}}

    \State Solve finite-horizon optimization problem using MPC: $\bigl(\mathbf{x}^\star, \mathbf{u}^\star\bigr)$
    \Comment{see~\eqref{eq:RL--MPC}}
    \State Compute objective of the predicted trajectory: $\hat J_{\text{RL--MPC}}\gets \hat J\bigl(\mathbf{x}^\star, \mathbf{u}^\star\bigr)$
    \Comment{see~\eqref{eq:objective-planning}}
    \State $u(t)\gets
    \begin{cases}
        u^\star(t), & \hat J_{\text{RL--MPC}}\ge \hat J_{\mathrm{RL}},\\
        \tilde u(t), & \text{otherwise},
    \end{cases}$
    \State Simulate next time step: $x(t+1)\gets f\!\left(x(t),u(t),d(t)\right)$ \Comment{see~\eqref{eq:greenlight-model}}
\EndFor
\end{algorithmic}
\end{algorithm}

\subsection{Performance metrics}\label{sec:performance-metrics}
Because this work combines and compares control methods using different optimization objectives, a discounted RL objective and an undiscounted finite-horizon MPC objective, closed-loop performance was evaluated using a common undiscounted metric.
This metric assumes that, under equal product prices and quality, one kilogram of tomatoes harvested today has the same economic value as one kilogram harvested later.

Closed-loop controller performance was quantified using three metrics computed over the entire simulation period.
The first metric is the cumulative objective $\mathcal J$, defined as the cumulative undiscounted sum of the reward~\eqref{eq:reward} over the controlled period $T$.
This metric quantifies how well each controller optimized the given economic control problem in~\eqref{eq:optimization_problem}.
The cumulative objective for a given controller is defined as:
\begin{equation}\label{eq:cumulative-objective}
    \mathcal J_{j} = \sum_{t=0}^{T-1} r(x(t), x(t+1), u(t))
\end{equation}
where $u(t)$ denotes the control input selected by controller $j\in\{\mathrm{MPC},\mathrm{RL},\mathrm{RL\text{--}MPC}\}$ at time step $t$, and the corresponding state sequence evolves according to the discrete-time model $f(\cdot)$ in~\eqref{eq:greenlight-model}.
The next two metrics decompose~\eqref{eq:cumulative-objective} into its economic component~\eqref{eq:economic-reward}, and penalty component~\eqref{eq:penalty-cost}. 
This yields the Economic Performance Indicator (EPI):
\begin{equation}\label{eq:cumulative-epi}
    \mathcal E_{j} = \sum_{t=0}^{T-1} r_{\mathrm{e}}(x(t), x(t+1), u(t)),
\end{equation}
and the cumulative penalty:
\begin{equation}\label{eq:cumulative-pen}
    \mathcal P_{j} = \sum_{t=0}^{T-1} r_{\mathrm{p}}(x(t+1)).
\end{equation}

\paragraph{Paired evaluation across weather trajectories.}
A paired evaluation was used to compare trajectory-selection RL--MPC against standalone MPC and RL under identical weather disturbances.
Specifically, for each weather trajectory  $\mathbf{d}$, the difference in a performance metric (i.e.,~\eqref{eq:cumulative-objective},~\eqref{eq:cumulative-epi}, or~\eqref{eq:cumulative-pen}) was computed and then averaged across trajectories.
For a finite weather set $\mathcal{D}$, the paired difference in cumulative objective is defined as:
\begin{equation}\label{eq:signed-difference}
  \Delta \mathcal{J} \;=\;
  \frac{1}{|\mathcal{D}|}
  \sum_{\mathbf{d}\in\mathcal{D}}
  \Bigl(
  \mathcal{J}_{\text{RL--MPC}}(\mathbf{d}) - \mathcal{J}_{j}(\mathbf{d})
  \Bigr).
\end{equation}
and analogously for $\Delta\mathcal{E}$ and $\Delta\mathcal{P}$ using~\eqref{eq:cumulative-epi} and~\eqref{eq:cumulative-pen}.

\paragraph{Comparing rollout and predicted trajectories of trajectory-selection RL--MPC.}
To analyze the selection mechanism in trajectory-selection RL--MPC, the finite-horizon planning objective~\eqref{eq:objective-planning} for the predicted solution and the guiding-policy rollout were compared:
\begin{equation}\label{eq:delta-objective}
    \Delta\hat J(t)
    =
    \hat J\bigl(\mathbf{x}^\star(t), \mathbf{u}^\star(t)\bigr)
    -
    \hat J\bigl(\tilde{\mathbf{x}}(t), \tilde{\mathbf{u}}(t)\bigr),
\end{equation}
where $\bigl(\mathbf{x}^\star(t), \mathbf{u}^\star(t)\bigr)$ denotes the predicted trajectories by solving~\eqref{eq:RL--MPC} and $\bigl(\tilde{\mathbf{x}}(t), \tilde{\mathbf{u}}(t)\bigr)$ denotes the RL policy rollout trajectories at time step $t$.

\begin{table}[t]
\centering
\caption{Weather disturbances, i.e., uncontrollable input $d$ to the tomato greenhouse model~\eqref{eq:greenlight-model}.}
\label{tab:gh_weather}
\begin{tabular}{@{}lcp{8.2cm}@{}}
\toprule
\textbf{Symbol} & \textbf{Unit} & \textbf{Description} \\
\midrule
$d_{
\mathrm{iGlob}}$ & $\mathrm{W\,m^{-2}}$ & Solar global radiation \\
$d_{\mathrm{T}}$ & $\mathrm{^\circ C}$ & Outdoor temperature \\
$d_{\ce{CO2}}$  & $\mathrm{ppm}$ &Outdoor CO$_2$ concentration \\
$d_{\mathrm{wind}}$ & $\mathrm{m\, s^{-1}}$ & Wind speed \\
$d_{\mathrm{RH}}$ & $\mathrm{\%}$ & Outdoor relative humidity \\
$d_{\mathrm{T_{sky}}}$ & $\mathrm{^\circ C}$ & Sky temperature \\
$d_{\mathrm{T_{soil}}}$ & $\mathrm{^\circ C}$ & Outdoor soil temperature \\
\bottomrule
\end{tabular}
\end{table}

\subsection{Weather datasets}\label{sec:weather-data}
Throughout this work, the dynamical model in~\eqref{eq:greenlight-model} was simulated using four finite datasets of historic weather trajectories.
Table~\ref{tab:gh_weather} lists the seven disturbance variables that the model takes as uncontrollable input.
Two datasets were used for policy training, and two were used for controller evaluation.

The specialist policy was trained on a single weather trajectory corresponding to June~1,~2023.
This specialist-policy training dataset is defined as:
\begin{equation}\label{eq:specialist-weather}
    \mathcal{D}^{\mathrm{train}}_{\mathrm{specialist}}\coloneq \left\{\mathbf{d}_{\mathrm{June\,1,\,2023}}\right\}
\end{equation}

The generalist RL policy was trained across multiple weather trajectories sampled from the period April~1 to June~1 over six years (2013--2018).
This generalist-policy training dataset is denoted by:
\begin{equation}\label{eq:generalist-trainingweather}
    \mathcal{D}^{\mathrm{train}}_{\mathrm{generalist}}
    \coloneq
    \left\{\mathbf{d}_{\mathrm{Apr\text{--}Jun},\,y}\; \middle|\; y\in\{2013,\dots,2018\}\right\},
\end{equation}
where each $\mathbf{d}$ represents an individual daily weather trajectory.

The specialist policy, standalone MPC controller, and trajectory-selection RL--MPC guided by the specialist policy were evaluated on the training weather trajectory: 
\[\mathcal{D}^{\mathrm{test}}_{\mathrm{specialist}}=\mathcal{D}^{\mathrm{train}}_{\mathrm{specialist}}.\]
The generalist policy, standalone MPC controller, and trajectory-selection RL--MPC guided by the generalist policy were evaluated on 25 individual dates, denoted by $\mathcal{D}^{\mathrm{test}}_{\mathrm{generalist}}$, selected from years held out during training (2011, 2012, 2019, 2020, and 2023).
Accordingly, the weather trajectories encountered during evaluation are disjoint from the generalist training set $\mathcal{D}^{\mathrm{train}}_{\mathrm{generalist}}\cap\mathcal{D}^{\mathrm{test}}_{\mathrm{generalist}}=\emptyset$.
The specific dates included in $\mathcal{D}^{\mathrm{test}}_{\mathrm{generalist}}$ are listed in Table~\ref{tab:generalist_test_dates}.

Weather data for 2011--2020 consists of real-world hourly measurements recorded at Schiphol Amsterdam, The Netherlands \citep{knmi_uurgeg_schiphol_2011_2020}.
These measurements were interpolated to a five-minute resolution using cubic interpolation.
Weather data from 2023 originates from the Autonomous Greenhouse Challenge (AGC) at Bleiswijk 2023, The Netherlands \citep{agc_2023}, and was already measured at five-minute intervals.

\begin{table}[t]
\centering
\caption{Evaluation dates for the generalist-policy dataset $\mathcal{D}^{\mathrm{test}}_{\mathrm{generalist}}$.
Day index is zero-based (day 0 = 1 Jan). Dates selected by day-of-year differ in leap years (2012, 2020).}
\label{tab:generalist_test_dates}
\begin{tabular}{lccccc}
\toprule
 & \multicolumn{5}{c}{\textbf{Day of year}} \\
\cmidrule(lr){2-6}
\textbf{Year} & \textbf{90} & \textbf{105} & \textbf{120} & \textbf{135} & \textbf{151} \\
\midrule
2011 & 01 Apr & 16 Apr & 01 May & 16 May & 01 Jun \\
2012 & 31 Mar & 15 Apr & 30 Apr & 15 May & 31 May \\
2019 & 01 Apr & 16 Apr & 01 May & 16 May & 01 Jun \\
2020 & 31 Mar & 15 Apr & 30 Apr & 15 May & 31 May \\
2023 & 01 Apr & 16 Apr & 01 May & 16 May & 01 Jun \\
\bottomrule
\end{tabular}
\end{table}

\subsection{Simulation experiments description}\label{sec:simulation-experiments}
The closed-loop performance of the proposed trajectory-selection RL--MPC framework was evaluated in two sets of simulation experiments, each guided by a different RL policy: the \emph{specialist policy}, and the \emph{generalist policy}.
To obtain numerical results within computationally feasible time, all closed-loop simulations spanned a fixed period of one day, resulting in $T=288$ time steps.
This evaluation period was shorter than a full greenhouse tomato production cycle.
Therefore, the results indicate whether trajectory-selection RL--MPC incorporates longer-term economic information into a short-horizon MPC optimization problem and improves closed-loop performance.
All the MPC optimization problems used a one-hour prediction horizon ($N=12$).
Table~\ref{tab:sim-settings} details all the simulation settings used in the finite-horizon optimization problems.

RL policies were trained offline using PPO, as described in Section~\ref{sec:rl}, i.e., by maximizing the expected discounted objective in~\eqref{eq:rl-objective}, with a single random seed.
The training and evaluation weather datasets are defined in Section~\ref{sec:weather-data}.
The coefficients in the economic reward in~\eqref{eq:economic-reward} and penalty function~\eqref{eq:penalty_factors} are listed in Table~\ref{tab:sim-settings}.
The indoor climate outputs ($y_{\mathrm{T, air}}$, $y_{\mathrm{CO_{2}, air}}$, $y_{\mathrm{RH, air}}$) were constrained with upper and lower bounds that reflect practical operational ranges and are listed in Table~\ref{tab:gh_state_bounds}.
The next two subsections detail the specialist policy and generalist policy experiment sets.

\begin{table}[t]
    \setlength{\tabcolsep}{4pt}   
    \footnotesize                 
    \caption{Economic reward, linear penalty function coefficients, simulation settings.}
    \label{tab:sim-settings}
    \begin{tabular*}{\columnwidth}{@{\extracolsep{\fill}} l l l p{0.42\columnwidth}}
    \toprule
    \textbf{Variable} & \textbf{Value} & \textbf{Unit} & \textbf{Description} \\
    \midrule
    $c_{\mathrm{CO_2}}$            & 0.30          & €/kg         & \ce{CO2} price \\
    $c_{\mathrm{heat}}$            & 0.09          & €/kWh       & Heating price \\
    $c_{\mathrm{elec}}$            & 0.20          & €/kWh       & Electricity price \\
    
    $c_{\mathrm{frt}}$            & 20.-      & €/kg\{DW\} & Dry matter fruit weight price\tablefootnote{Corresponding to a fresh weight price of 1.2 €/kg} \\
    $\lambda_{\mathrm{CO_2}}$      & $5\times10^{-5}$& --                   & Coefficient for CO\textsubscript{2} violations \\
    $\lambda_{\mathrm{T}}$      & $5\times10^{-3}$& --                   & Coefficient for temperature violations \\
    $\lambda_{\mathrm{RH}}$        & $7\times10^{-4}$& --                   & Coefficient for relative humidity violations \\
        $\Delta t$ & $300$ & $\mathrm{s}$ & Control interval\\
    $N$ & $12$ &  -- & Prediction horizon (time steps)\\
    $T$ & $288$ & -- & Simulation length (time steps)\\
    \bottomrule
    \end{tabular*}
\end{table}

\begin{table}[t]
    \caption{Upper and lower boundaries for the three indoor climate variables.}
    \label{tab:gh_state_bounds}
    \begin{tabular}{lcccl}
        \toprule
        \textbf{Symbol} & \textbf{Lower bound} ($y^{\min}_{i}$)& \textbf{Upper bound} ($y^{\max}_{i}$) & \textbf{Unit} & \textbf{Description} \\
        \midrule
        $y_{
        \mathrm{T, air}}$ & $10$ & $25$ & $\mathrm{^\circ C}$ & Greenhouse air temperature \\
        $y_{\mathrm{CO_2, air}}$ & $400$ & $1600$ & $\mathrm{ppm}$ & Carbon dioxide concentration \\
        $y_{\mathrm{RH, air}}$ & $0$ & $90$ & $\mathrm{\%}$ & Relative humidity \\
        \bottomrule
    \end{tabular}
\end{table}
\subsubsection*{Trajectory-selection RL--MPC guided by the specialist policy}\label{sec:sim-exp-specialist}
Three simulation experiments evaluated the performance and the behavior of the trajectory-selection RL--MPC framework when guided by the specialist RL policy.
This policy was trained on a single weather trajectory $\mathcal{D}^{\mathrm{train}}_{\mathrm{specialist}}$ (defined in~\eqref{eq:specialist-weather}), yielding a high-performing specialist RL policy.
All controllers were evaluated on that same weather trajectory $\mathcal{D}^{\mathrm{test}}_{\mathrm{specialist}}=\mathcal{D}^{\mathrm{train}}_{\mathrm{specialist}}$.
\begin{enumerate}
    \item First, the three resulting controllers (RL, MPC, and trajectory-selection RL--MPC guided by the specialist policy) were evaluated using the closed-loop performance metric: the cumulative objective~\eqref{eq:cumulative-objective}, the EPI~\eqref{eq:cumulative-epi}, and the cumulative penalty~\eqref{eq:cumulative-pen}.
    This experiment assessed whether trajectory-selection RL--MPC can incorporate longer-term economic information from a high-performing policy into a short-horizon MPC controller.
    A breakdown of the underlying costs and penalties provided insight into which strengths of the underlying methods contributed to controller performance.
    \item Second, the selection mechanism of trajectory-selection RL--MPC was analyzed by evaluating $\Delta \hat J$ in~\eqref{eq:delta-objective} over the simulation day.
    This analysis should show whether the input $u(t)$ from the RL rollout trajectory or from the MPC-predicted trajectory was selected and applied at each time step.
    It is hypothesized that the RL rollout trajectory would be selected more often during daytime, when solar radiation requires balancing current operating costs against future fruit yield.
    \item Third, an ablation study removed three elements of the trajectory-selection RL--MPC framework (Algorithm~\ref{alg:rlmpc_closedloop}) one at a time: the terminal cost $\hat{V}_{f}$, the terminal region constraint $\mathbb{X}_{f}$, and the selection mechanism.
    Each ablated version of the framework was evaluated by the cumulative objective, EPI, and cumulative penalty on $\mathcal{D}^{\mathrm{test}}_{\mathrm{specialist}}$.
    When the selection mechanism was removed, the first input was always taken from the MPC-predicted trajectory $\mathbf{u}^\star$.
    This experiment assessed which components of the framework are necessary for the observed performance gains.
\end{enumerate}

\subsubsection*{Trajectory-selection RL--MPC guided by the generalist policy}
The second set of simulation experiments evaluated trajectory-selection RL--MPC when guided by the generalist RL policy on weather trajectories held out during policy training.
This policy was trained by sampling weather trajectories uniformly at random from $\mathcal{D}^{\mathrm{train}}_{\mathrm{generalist}}$ (described in Section~\ref{sec:weather-data}).
The three controllers (RL, MPC, and trajectory-selection RL--MPC guided by the generalist policy) were evaluated on 25 individual test dates, denoted by $\mathcal{D}^{\mathrm{test}}_{\mathrm{generalist}}$.
The specific dates are listed in Table~\ref{tab:generalist_test_dates}.
Controller performance was compared using the paired evaluation metric in~\eqref{eq:signed-difference} for the cumulative objective, EPI, and cumulative penalty.
This simulation experiment aimed to assess whether trajectory-selection RL--MPC can improve closed-loop performance relative to the guiding RL policy and MPC on weather trajectories not encountered during policy training.
Although training RL policies on broader weather datasets can improve generalization, it does not prevent degraded performance under unseen weather conditions.
Therefore, it is hypothesized that trajectory-selection RL--MPC could further improve over the guiding RL policy through finite-horizon MPC optimization using short-term actual weather predictions.
In addition, the experiment assessed when the proposed framework benefits most from the guiding RL policy or from short-horizon MPC optimization.
The relationships between daily mean outdoor temperature, global radiation, and relative humidity, and both the paired performance difference $\Delta \mathcal{J}$ and the trajectory-selection mechanism, were analyzed.
This analysis assessed whether unseen weather conditions outside or near the boundary of the training distribution were associated with larger performance differences between trajectory-selection RL--MPC and the guiding RL policy, and whether the selection mechanism relied more strongly on the RL rollout or the MPC-predicted trajectory under specific weather conditions.

\section{Results}\label{sec:results}
\begin{figure}
    \centering
    \includegraphics[width=1\linewidth]{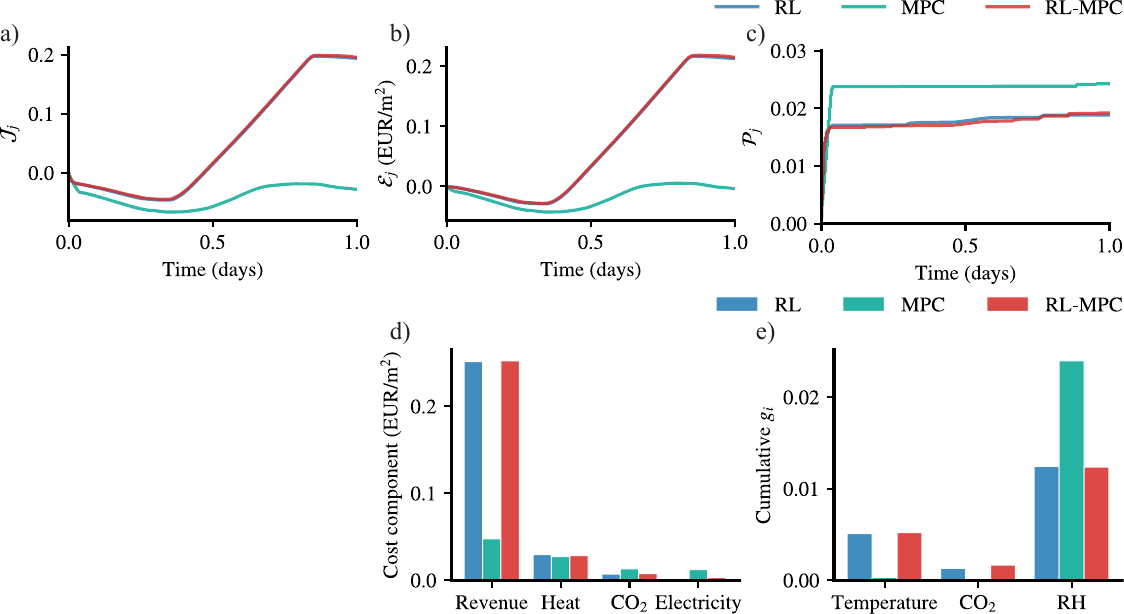}
    \caption{\textbf{Comparison of closed-loop controller performance on the specialist-policy training trajectory $\mathcal{D}^{\mathrm{train}}_{\mathrm{specialist}}$.}
    This figure visualizes the cumulative sums of three performance metrics and their individual components for MPC, RL, and trajectory-selection RL--MPC, evaluated on a single weather trajectory.
    Note that the cumulative sums of RL and RL--MPC in a), b), and c) are nearly identical and, therefore, have overlapping trajectories.
    The weather trajectory corresponds to the trajectory used to train the specialist RL policy, i.e., $\mathcal{D}^{\mathrm{test}}_{\mathrm{specialist}}=\mathcal{D}^{\mathrm{train}}_{\mathrm{specialist}}$.
    \textbf{a)} Cumulative objective ($\mathcal{J}$) as defined in~\eqref{eq:cumulative-objective}. 
    \textbf{b)} EPI ($\mathcal{E}$) as defined in~\eqref{eq:cumulative-epi}. 
    \textbf{c)} Cumulative penalty ($\mathcal{P}$) as defined in~\eqref{eq:cumulative-pen}.
    \textbf{d)} Cumulative components of the economic reward as defined in~\eqref{eq:economic-reward}.
    \textbf{e)} Cumulative penalty of the individual penalty components $g_i$, as defined in~\eqref{eq:penalty_factors}.
    The objective and EPI are maximized, whereas penalty values are minimized.
    For the economic components, higher revenue and lower operating costs indicate better performance.
    }
    \label{fig:specialist-cost-june-1}
\end{figure}

\subsection{Trajectory-selection RL--MPC guided by the specialist policy}\label{sec:res-specialist}
The results for the specialist RL policy, MPC, and trajectory-selection RL--MPC guided by the specialist policy evaluated on the specialist policy evaluation set $\mathcal{D}^{\mathrm{test}}_{\mathrm{specialist}}=\mathcal{D}^{\mathrm{train}}_{\mathrm{specialist}}$ are presented in Figure~\ref{fig:specialist-cost-june-1}.
Trajectory-selection RL--MPC, with a one-hour prediction horizon, matched the specialist policy's cumulative performance $\mathcal J$.
The cumulative objective, EPI, and cumulative penalty exhibited nearly identical cost trajectories over the optimized day, as shown in Figure~\ref{fig:specialist-cost-june-1}.
The trajectory-selection RL--MPC controller provided only a marginal improvement in the cumulative objective $\mathcal J$ ($0.196$ vs $0.194$) compared to the specialist RL policy.
This gain was driven by a relatively small improvement in the EPI ($0.215$ vs $0.213$), while the cumulative penalty was identical ($0.019$) for this specific weather trajectory.
The plots of the individual components of the economic reward and the penalty function confirm that trajectory-selection RL--MPC performed nearly identically to the specialist policy.
With the same prediction horizon, the trajectory-selection RL--MPC controller outperformed the MPC controller by a relatively large margin ($\Delta\mathcal J=0.224$).
This improvement was driven by a higher EPI ($0.196$ vs $-0.004$).
The higher EPI was mainly due to increased revenue under trajectory-selection RL--MPC compared with MPC ($+0.20$), while similar amounts of heating and \ce{CO2} injection were applied.
In addition, the MPC controller incurred a higher cumulative penalty than the trajectory-selection RL--MPC controller ($+0.005$).
This increase was primarily due to more severe violations of the relative humidity boundary at the start of the day.

\begin{figure}
\centering
\includegraphics[width=0.5\linewidth]{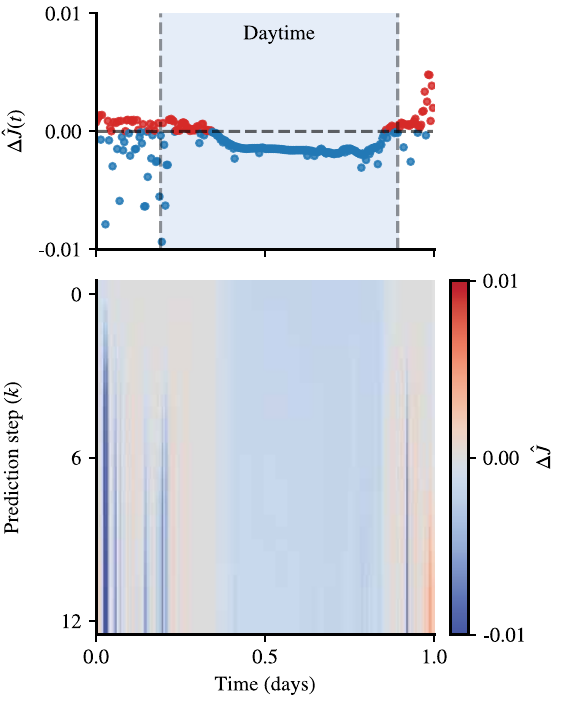}
    \caption{\textbf{Evaluating the finite-horizon planning objective of the policy rollout and the MPC-predicted trajectory for $\mathcal{D}^{\mathrm{train}}_{\mathrm{specialist}}$.}
    \textbf{Top plot:} Planning-objective difference $\Delta\hat J(t)$~\eqref{eq:delta-objective} between the rollout generated by $\piRL$ and the MPC-predicted solution over the prediction horizon $N$ at time step $t$, plotted against time of day.
    Red (positive) values indicate that the MPC-predicted trajectory (with terminal ingredients) achieved a higher planning objective, whereas blue (negative) values indicate that the policy rollout achieved a higher planning objective.
    \textbf{Bottom plot:} Heatmap of the planning-objective difference between the rollout and predicted trajectories evaluated up to prediction step $k$, plotted against time of day.
    Again, red indicates that the MPC-optimization achieved a higher planning objective at prediction step $k$, whereas blue indicates that the RL rollout achieved a higher planning objective.
    }F
\label{fig:decision-heatmap}
\end{figure}

At each time step, the first control input from the trajectory (rollout or predicted) with the best finite-horizon planning objective~\eqref{eq:objective-planning} was selected and applied to the system.
Evaluating the difference in the planning objective $\Delta \hat J$~\eqref{eq:delta-objective} between the policy rollout and the MPC-predicted trajectory (with terminal ingredients) shows that the specialist policy found better solutions during a consistent period of the daytime, see Figure~\ref{fig:decision-heatmap}.
This period starts a few hours after sunrise, i.e., the start of daytime, and ends around sunset.
This period also corresponds to the period during which both RL and trajectory-selection RL--MPC improved their cumulative objective the most compared to MPC in Figure~\ref{fig:specialist-cost-june-1}(a).
The heatmap confirms this behavior, as the policy rollout had a higher objective at each prediction step ($k$) during the daytime.
Conversely, at the end of the simulated day, the MPC-predicted trajectories yielded better solutions, indicated by a higher finite-horizon planning objective $\hat J$.
At the start of the simulated day, trajectory-selection RL--MPC frequently switched between selecting the input from the rollout or the predicted trajectory.

The ablation study revealed that removing the selection mechanism had the largest impact on controller performance, as the cumulative objective dropped by $12.2\%$ from ($0.196$) to ($0.172$).
This performance degradation is mainly due to the cumulative penalty, which increased by $95.2\%$ from ($0.021$) to ($0.041$).
This increase in cumulative penalty occurred primarily at the beginning and the end of the simulated day, when lower temperature and higher relative humidity drove the indoor climate variables closer to their constraint bounds, making penalties more likely.
On the other hand, the EPI did not drop significantly when discarding the input selection mechanism.
Removing the terminal cost and the terminal region constraint had minimal effect on the trajectory-selection RL--MPC controller's performance when guided by the specialist policy.

\begin{figure}
    \centering
    \includegraphics[width=\linewidth]{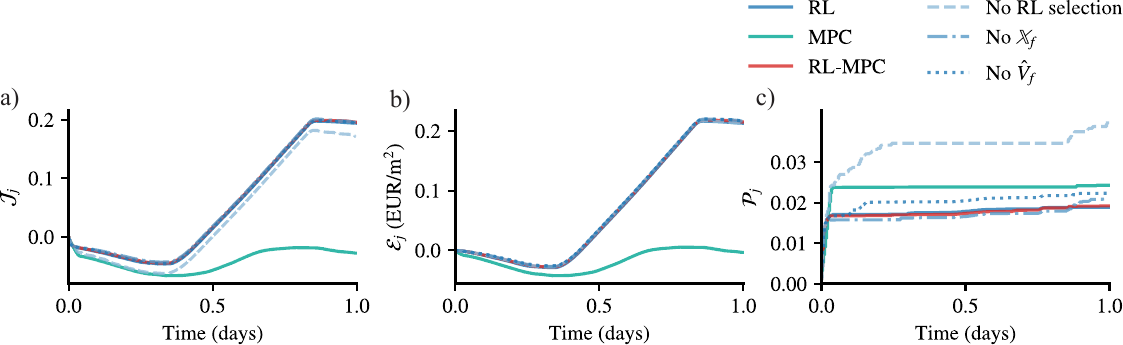}
    \caption{\textbf{Ablation study for trajectory-selection RL--MPC.} 
    These graphs show the cumulative performance metrics of MPC, RL, trajectory-selection RL--MPC, and ablated variants of the proposed framework obtained by removing individual components from the framework as defined in Algorithm~\ref{alg:rlmpc_closedloop}.
    The resulting controllers were evaluated on the specialist-policy training trajectory, i.e., $\mathcal{D}^{\mathrm{test}}_{\mathrm{specialist}}=\mathcal{D}^{\mathrm{train}}_{\mathrm{specialist}}$.
    \textbf{a)} Cumulative objective as defined in~\eqref{eq:cumulative-objective}. 
    \textbf{b)} EPI ($\mathcal{E}$) as defined in~\eqref{eq:cumulative-epi}. 
    \textbf{c)} Cumulative penalty ($\mathcal{P}$) as defined~\eqref{eq:cumulative-pen}. 
    Here, $\mathbb{X}_f$ denotes the terminal region constraint and $\hat{V}_{f}$ denotes the terminal cost.
    }
    \label{fig:ablation}
\end{figure}

\subsection{Trajectory-selection RL--MPC guided by the generalist policy}\label{sec:res-generalist}
\begin{figure}
    \centering
    \includegraphics[width=\linewidth]{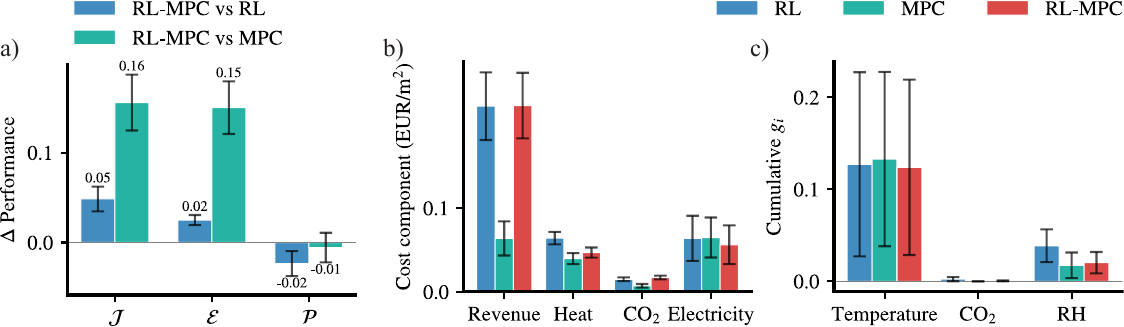}
    \caption{\textbf{Comparison of closed-loop performance of RL, MPC, and trajectory-selection RL--MPC guided by the generalist RL policy.} 
    Bars show mean values across 25 simulated dates ($\mathcal{D}^{\mathrm{test}}_{\mathrm{generalist}}$) with $95\%$ confidence intervals. 
    \textbf{a)} Mean signed differences between RL--MPC and each standalone controller, computed using~\eqref{eq:signed-difference}.
    \textbf{a)} Mean signed difference between RL--MPC and each standalone controller, computed using~\eqref{eq:signed-difference}.
    Positive values of the objective $\mathcal J$ and EPI indicate that RL--MPC outperforms the corresponding standalone controllers.
    For the penalty metric, negative values indicate improvement for RL--MPC over the corresponding controller (i.e., fewer constraint violations).
    \textbf{b)} Mean cumulative values of the individual components of the economic reward function defined in~\eqref{eq:economic-reward}.
    \textbf{c)} Mean cumulative penalty of the individual penalty components $g_i$, as defined in~\eqref{eq:penalty_factors}.}
    \label{fig:generalist-bars}
\end{figure}

Figure~\ref{fig:generalist-bars} summarizes the simulation results on the evaluation set $\mathcal{D}^{\mathrm{test}}_{\mathrm{generalist}}$ for the generalist RL policy, MPC, and trajectory-selection RL--MPC guided by the generalist policy.
Across all 25 evaluated dates, trajectory-selection RL--MPC outperformed both standalone controllers.
On average, trajectory-selection RL--MPC increased performance by $0.048\pm0.013$ ($+53.9\%$) relative to the generalist RL policy.
This improvement was driven by a higher EPI ($+0.025\,\mathrm{EU/m^2}$) and a lower cumulative penalty ($-0.024$).
An analysis of the economic reward components shows that the trajectory-selection RL--MPC controller matched the revenue achieved by the generalist RL policy ($0.222\,\mathrm{EU/m^2}$), while reducing heating costs to $0.047\,\mathrm{EU/m^2}$ ($-26.8\%$) and lamp electricity costs to $0.056\,\mathrm{EU/m^2}$ ($-12.0\%$).
Decomposing the linear penalty further indicates that trajectory-selection RL--MPC mainly benefits from fewer relative humidity violations, with a cumulative penalty of $0.020\pm0.012$ ($-47.7\%$).

Compared to standalone MPC, trajectory-selection RL--MPC increased the cumulative objective by a significant margin of $0.156\pm0.031$ ($+79.9\%$).
This gain was mainly due to an increase in EPI of $0.150\,\mathrm{EU/m^2}$.
Examining the cost components suggests that the discrepancy in EPI between the two controllers is explained by the low revenue achieved by MPC ($0.063\,\mathrm{EU/m^2}$). 
Trajectory-selection RL--MPC also spent less on electricity for the lamps ($-13.8\%$).
The difference in cumulative penalty between the trajectory-selection approach and MPC is relatively small and not significant $(-0.006\pm0.016)$.
Also note the large confidence intervals for the cumulative penalty on temperature violations for all three controllers.

\begin{figure}
    \centering
    \includegraphics[width=\linewidth]{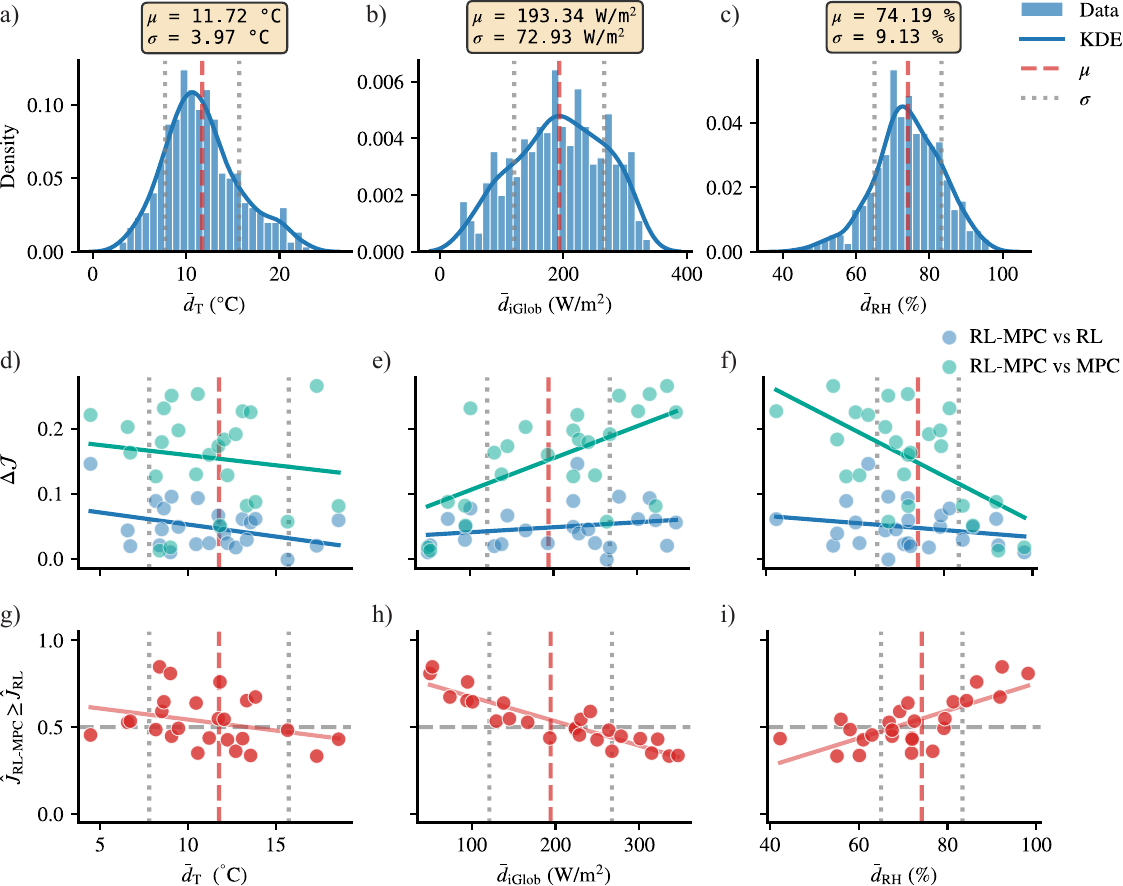}
    \caption{\textbf{Relation between daily mean weather variables, trajectory-selection RL--MPC performance gains, and the input selection mechanism.} 
    \textbf{Top row:} Empirical distributions and kernel density estimates (KDEs) of the daily mean weather variables in the training set $\mathcal{D}^{\mathrm{train}}_{\mathrm{generalist}}$: a) outdoor temperature $\bar d_{\mathrm{T}}$, b) global radiation $\bar d_{\mathrm{iGlob}}$, c) relative humidity \textbf{$\bar d_{\mathrm{RH}}$}. 
    Daily mean was computed for all $366$ days in this set. 
    The red dashed line indicates the mean $\mu$, and the gray dotted lines represent $\mu\pm\sigma$.
    \textbf{Middle row:} Closed-loop improvement $\Delta\mathcal{J}$~\eqref{eq:signed-difference} of RL--MPC over the generalist RL policy and MPC on individual days in the test set $\mathcal{D}^{\mathrm{test}}_\mathrm{generalist}$.
    The x-axis shows the daily mean weather variables for each evaluated day: d) outdoor temperature $\bar d_{\mathrm{T}}$, e) global radiation $\bar d_{\mathrm{iGlob}}$, f) relative humidity \textbf{$\bar d_{\mathrm{RH}}$}.
    Solid lines indicate linear trend fits.    \textbf{Bottom row:} Fraction of time the predicted trajectory attained a higher finite-horizon planning objective $\hat J$~\eqref{eq:objective-planning} than the RL rollout trajectory.
    The x-axis shows the daily mean weather variables for each evaluated day: g) outdoor temperature $\bar d_{\mathrm{T}}$, h) global radiation $\bar d_{\mathrm{iGlob}}$, i) relative humidity \textbf{$\bar d_{\mathrm{RH}}$}.
    Solid lines indicate linear trend fits.
    }
    \label{fig:generalist-mean-weather}
\end{figure}

Figure~\ref{fig:generalist-mean-weather}(d--f) plots the performance improvements over the two standalone controllers against the the daily mean outdoor temperature, global radiation, and relative humidity ($\bar d_{\mathrm{T}},\,\bar d_{\mathrm{iGlob}},\,\bar d_{\mathrm{RH}}$) for each day in $\mathcal{D}^{\mathrm{test}}_\mathrm{specialist}$. 
One can observe that trajectory-selection RL--MPC outperforms both standalone controllers on all evaluation days, as indicated by $\Delta \mathcal J>0$.
In particular, the performance difference between trajectory-selection RL--MPC and MPC shows a strong positive linear relation with global radiation.
Higher global radiation corresponds to larger improvements of trajectory-selection RL--MPC over MPC.
In contrast, the performance difference between trajectory-selection RL--MPC and MPC shows a negative linear relation with the daily mean relative humidity, indicating that the performance gap between the two methods narrows under higher relative humidity.
The selection mechanism shows consistent trends with these variables, as shown in Figure~\ref{fig:generalist-mean-weather}(g--i): higher daily average global radiation increased the fraction of time the guiding policy rollout yielded the better solution to~\eqref{eq:objective-planning}, leading to selection of its first control input.
Conversely, higher daily mean relative humidity increased the fraction of time the predicted trajectory yielded a better solution than the guiding policy rollout.
Temperature did not seem to have a significant effect on the relative performance or the selection mechanism.
The improvement of trajectory-selection RL--MPC over the generalist RL policy appeared less sensitive to the daily mean weather variables, even on test days whose weather statistics fall outside the training distribution, see Figure~\ref{fig:generalist-mean-weather}(a--c).

\section{Discussion}\label{sec:discussion}
\subsection{Trajectory-selection RL--MPC guided by the specialist policy}\label{sec:discussion-specialist}
The closed-loop performance of trajectory-selection RL--MPC with a one-hour prediction horizon matched the performance of the specialist policy when evaluated on the weather trajectory the policy was trained on~\eqref{eq:specialist-weather}.
This result suggests that the proposed framework can incorporate longer-term economic information, encoded by the guiding RL policy, into a short-horizon MPC optimization problem.
As expected, trajectory-selection RL--MPC outperformed MPC with the same horizon length by a large margin $\Delta \mathcal J=0.224$.
Figure~\ref{fig:specialist-cost-june-1} shows that this performance improvement was mainly driven by increased EPI. 
This shows that the one-hour prediction horizon is too short for MPC to balance current operating costs against the delayed benefits of crop yield.
This inability of the myopic MPC was reflected in the low revenue despite only minor differences in heating, \ce{CO2}, and electricity costs compared with RL and trajectory-selection RL--MPC.
This observation suggests that MPC prioritized short-term objectives, such as minimizing constraint violations and operating costs, potentially by using the thermal and blackout screens at the expense of crop growth.
Such control strategies limit PAR entering the greenhouse, thereby inhibiting crop yield.
In contrast, trajectory-selection RL--MPC achieved the same revenue as the specialist RL policy despite using the same one-hour MPC prediction horizon.

These results demonstrate a general limitation of short-horizon MPC control for greenhouse tomato production systems.
A one-hour prediction horizon for receding horizon approaches, such as MPC, can lead to myopic behavior in greenhouse tomato production systems.
Such horizon lengths are typically too short to balance immediate operating costs against the delayed benefit of crop yield.
Ideally, MPC would use prediction horizons spanning the entire production period \citep{VANHENTEN2009timescale, vanstratenOptimalGreenhouseCultivation2010, kuijpersFruitDevelopmentModelling2021}.
However, as greenhouse fruit production models become more complex, and as uncertainties in weather forecasts, product prices, and crop dynamics are incorporated, optimizing over such long horizons becomes computationally intractable.
This problem setting motivates control methods that incorporate longer-term economic information into finite-horizon control methods.
In this work, trajectory-selection RL--MPC improved closed-loop performance over MPC with the same prediction horizon through terminal region constraints, terminal costs, and trajectory selection.

All simulations in this work spanned a single day.
In practice, greenhouse tomato production cycles have a duration of weeks \citep{dekoningDevelopmentDryMatter1994, adamsEffectTemperatureGrowth2001}.
For instance, fruit growth is influenced by daily mean crop temperature, and controlling the system for multiple consecutive days may require longer prediction horizons.
Moreover, the employed tomato yield model does not account for fruit development stages: all carbohydrates partitioned to fruits, minus those used for maintenance respiration, were counted as yield.
Ignoring the temporal delay of fruit development may overestimate short-term yield \citep{katzinGreenLightOpenSource2020}.
Incorporating fruit development into the model dynamics causes the controller's decisions to affect fruit yield for much longer periods.
The presented results show that trajectory-selection RL--MPC can incorporate longer-term information into a short-horizon optimization problem over a closed-loop simulation period that is significantly longer than the MPC prediction horizon.
Future research should evaluate if the guiding policy can capture performance across full production cycles with explicit fruit development dynamics, and transfer this longer economic information effectively to MPC optimization.
In addition, future work should compare the proposed framework to timescale decomposition approaches \citep{VANHENTEN2009timescale, tapEconomicsbasedOptimalControl2000}.

The benefit of the guiding policy is further demonstrated by comparing the finite-horizon planning objective $\hat J$ in $\eqref{eq:objective-planning}$ between the trajectories from the policy rollout and the trajectory-selection RL--MPC optimization.
For a consistent period during the daytime, the policy rollout attained a higher planning objective than the predicted trajectory; see Figure~\ref{fig:decision-heatmap}.
In this period, global radiation induced crop growth, and trajectory-selection RL--MPC therefore selects the control input from the policy rollout since the policy benefits from its relatively long effective horizon (approximately 16 hours; Section~\ref{sec:rl}) compared with the one-hour optimization horizon.
This interval also coincides with the steep rise in the cumulative objective and EPI observed in Figure~\ref{fig:specialist-cost-june-1}.
Near the end of the simulated day, the MPC-predicted trajectories yielded better solutions, as indicated by a positive $\Delta \hat J$.
After sunset, the indoor temperature dropped, and relative humidity rose, increasing the importance of constraint handling, a known strength of MPC.

The decreased performance of standalone MPC regarding the constraint violations, shown in Figure~\ref{fig:specialist-cost-june-1}, and the relatively frequent selection of the RL rollout trajectory, shown in Figure~\ref{fig:decision-heatmap}, may be partly explained by the relaxed optimization tolerances used in the finite-horizon MPC problem, as mentioned in Section~\ref{sec:mpc}.
With these settings, the nonlinear optimizer may terminate before converging to a locally optimal solution, with a worse finite-horizon planning objective value than the RL rollout.
As a result, MPC-predicted trajectories may have been rejected by the trajectory-selection mechanism, even though stricter tolerances might have produced improved MPC solutions.
The relaxed tolerances were used to keep the closed-loop simulations computationally feasible.
Imposing stricter tolerances on the optimization problem could improve the performance of MPC-predicted trajectories, possibly increasing the benefit of trajectory-selection RL--MPC.
However, this would increase the controller's computational demand and does not eliminate the limitation imposed by the short prediction horizon.

This work trained the RL policy using a discounted objective by introducing the discount factor $\gamma=0.995$, as described in Section~\ref{sec:rl}.
For the considered economic greenhouse control problem, the true performance metric is undiscounted, since economic returns (e.g., tomato yield) obtained later in the production period are not worth less under equal prices and product quality.
Ideally, this would correspond to RL using $\gamma=1$.
However, training RL policies using deep neural networks with an undiscounted objective introduces high variance in gradient estimates, which can make policy optimization unstable.
To limit this variance and improve the convergence of RL policies, discounting was applied to compromise between longer-term economic returns and maintaining stable RL training.
Since short-horizon MPC lacks knowledge of the future return, RL trained with a discounted objective still encodes longer-term economic information, while distorting the value of future returns and therefore changes the objective for which RL was optimized.
Consequently, evaluating the RL policy, the MPC controller, and trajectory-selection RL--MPC using common undiscounted closed-loop metrics introduces an objective mismatch for RL
Similarly, this mismatch may affect the selection mechanism, which compares the MPC-predicted trajectory and the RL rollout using the undiscounted finite-horizon planning objective.
Therefore, the mechanism asserts the improvement under the undiscounted economic objective, not under the RL training objective.
This comparison may favor MPC-predicted trajectories.
Further research should evaluate how the value of the discount factor affects the learned policy, the trajectory-selection decisions, and the resulting closed-loop performance.

The ablation study highlights the role of the selection mechanism in achieving specialist-policy performance on the training trajectory.
Figure~\ref{fig:ablation} shows that removing this element increased the cumulative penalty, particularly near the start of the simulated day.
Figure~\ref{fig:decision-heatmap} further demonstrates the dependence of trajectory-selection RL--MPC on the policy rollout; at multiple time steps $t$ until sunrise, the rollout attained a substantially higher finite-horizon planning objective value than the predicted trajectory.
In this tested variant of trajectory-selection RL--MPC, the optimization using terminal ingredients appeared to be insufficient to maintain a low cumulative penalty and can even perform worse than standalone MPC in terms of cumulative penalty.
This behavior may be caused by the quasi-Newton approximation of the Hessian and the relatively loose tolerance afforded to the optimization routine, which were selected to ensure computational tractability but did not ensure that the optimization algorithms return local optima.
All three ablated versions performed identically on EPI, indicating that long-term economic guidance is preserved whenever the guiding policy is used, either through the selection mechanism or the terminal ingredients.

Overall, these results indicate that, on the specialist policy training trajectory, the performance of trajectory-selection RL--MPC is mainly driven by its selection mechanism.
Because the guiding policy was trained extensively on this specific weather trajectory, the RL policy rollout was already high-performing, so selecting the rollout trajectory provided most of the gains in EPI and cumulative penalty.
As a result, the terminal ingredients contributed less to controller performance in this setting.
However, under unseen weather conditions, the policy rollout trajectory is expected to be selected less often, and the contribution of the terminal ingredients may become more important.
This is discussed further for trajectory-selection RL--MPC guided by the generalist policy in Section~\ref{sec:discussion-generalist}.

Finally, all controllers were evaluated under the assumptions of full state availability and perfect knowledge of the actual future weather.
In practice, these measurements and forecasts are subject to error; therefore, the reported controller performance should be interpreted as an upper bound.
In addition, a single random seed is used to train the guiding policy.
Previous work on integrating RL and MPC indicates that the RL policy varies across random seeds, but that incorporating terminal ingredients in the trajectory-selection RL--MPC optimization framework can mitigate this variance \citep{vanlaatumStochasticModelPredictive2026}.

\subsection{Trajectory-selection RL--MPC guided by the generalist policy}\label{sec:discussion-generalist}
The simulation experiments presented in Section~\ref{sec:res-generalist} evaluated controller performance on an evaluation set $\mathcal{D}^{\mathrm{test}}_{\mathrm{generalist}}$ that differed from the generalist policy’s training set~\eqref{eq:generalist-trainingweather}.
On average, trajectory-selection RL--MPC improved the cumulative objective $\Delta \mathcal J$~\eqref{eq:signed-difference} by $0.048\pm0.013$ relative to the generalist RL policy and by $0.156\pm0.031$ relative to MPC.
The improvement over the generalist policy was driven by both a higher EPI and a lower cumulative penalty. 
The economic reward decomposition in Figure~\ref{fig:generalist-bars}(b) shows that trajectory-selection RL--MPC matched the revenue of the guiding policy while reducing heating and electricity costs.
In addition, trajectory-selection RL--MPC incurred fewer constraint violations, primarily due to a reduced penalty for relative humidity.
Interestingly, trajectory-selection RL--MPC also improved economic performance relative to the guiding policy by reducing operating costs.
These observations suggest that the receding-horizon optimization in trajectory-selection RL--MPC, with terminal ingredients, was able to refine the guiding policy rollout trajectory, thereby reducing operating costs and the cumulative penalty while maintaining revenue.
Compared with MPC, the increased cumulative objective was induced by a higher revenue, consistent with the specialist-policy results discussed in Section~\ref{sec:discussion-specialist}.
The strong dependence of temperature violations on external weather conditions led to relatively large confidence intervals and cumulative penalties for all three controllers, as shown in Figure~\ref{fig:generalist-bars}(c).
As outdoor temperature and global radiation increase, the controllers cannot sufficiently reduce indoor temperature, since the modeled greenhouse system lacks an active cooling system.

Analyzing the difference in cumulative objective further showed that trajectory-selection RL--MPC outperformed the standalone controllers on all 25 simulated days, see Figure~\ref{fig:generalist-mean-weather}. 
The performance difference between trajectory-selection RL--MPC and MPC exhibited a positive linear relation with daily mean global radiation.
Higher radiation induced crop growth, which increased the benefit of longer horizons and improved EPI performance, thereby amplifying the advantage over the myopic MPC controller.
In contrast, this performance difference showed a negative linear relation with daily mean outdoor relative humidity.
Elevated outdoor relative humidity increased the importance of maintaining the indoor relative-humidity constraint and reduced the emphasis on EPI, resulting in smaller performance differences between the two methods.
This effect is also illustrated in the selection mechanism in Figure~\ref{fig:generalist-mean-weather}.
Higher global radiation increased the fraction of time the policy rollout was selected, whereas higher daily mean humidity increased the fraction of time the predicted trajectory was selected.

The ablation study, discussed in the previous section, showed that, on the specialist-policy training trajectory, the selection mechanism was the primary contributor to controller performance.
However, under unseen weather conditions, the predicted MPC trajectory was more likely to be selected, as illustrated in Figure~\ref{fig:generalist-mean-weather}(g-i).
This observation suggests that, in such settings, the terminal ingredients play a more important role in controller performance than in the specialist case.
Previous work has also shown that policy-guided MPC with terminal constraints and terminal costs improves controller performance \citep{MSAAD2025449, vanlaatumStochasticModelPredictive2026}.
Hence, these findings support the need for terminal ingredients, in addition to the selection mechanism, to improve closed-loop controller performance under unseen weather conditions.

Overall, these results suggest that when the test set differs from the training set, trajectory-selection RL--MPC effectively leverages the complementary strengths of its underlying methods.
The controller retained the long-term economic performance of the guiding policy through high revenue due to crop growth, while the MPC improved short-term constraint handling and decreased operating costs.
The trajectory-selection RL--MPC controller was evaluated on weather trajectories not included in the RL policy's training set.
When the RL policy was trained on a broader weather dataset and faced unseen weather conditions at evaluation, its closed-loop performance may still decrease.
Trajectory-selection RL--MPC's receding-horizon approach can refine control input trajectories, reducing operating costs and constraint violations relative to the RL policy.
This mechanism is not limited to weather variability and may extend to other external disturbances, such as changes in crop prices and resource costs.
However, trajectory-selection RL--MPC does require a sufficiently accurate prediction model.
Modeling errors or uncertainty in the state dynamics can reduce, or even reverse, the benefit of the proposed framework.
When such uncertainties affect controller performance, uncertainty-aware approaches that combine RL with scenario-based stochastic MPC (SMPC) offer a feasible alternative.
For instance, RL-SMPC learns an RL policy in a stochastic setting and uses it to construct terminal ingredients and provide nonlinear feedback within scenario-based SMPC under parametric uncertainty \citep{vanlaatumStochasticModelPredictive2026}.

\section{Conclusion}\label{sec:conclusion}
The main contribution of this paper is \textit{trajectory-selection RL--MPC}, a framework designed to improve short-horizon MPC for greenhouse fruit production systems.
It was shown that the delayed economic benefits of fruit yield can be incorporated into a short-horizon MPC optimization problem by incorporating a terminal region constraint and terminal cost, derived from an RL rollout trajectory, and an online trajectory-selection mechanism.
In addition, the finite-horizon MPC optimization refined the control input trajectories using short-term weather disturbances, improving performance relative to the guiding policy by reducing constraint violations and operating costs when evaluated on weather trajectories not encountered during policy training.
The framework was applied on GreenLight, a large-scale, nonlinear greenhouse tomato model that exhibits stiff dynamics, for which optimization with long prediction horizons is computationally demanding.

The simulation results demonstrate that trajectory-selection RL--MPC with a one-hour prediction horizon matched the closed-loop performance of a specialist policy trained on a single weather trajectory.
Under this setting, trajectory-selection RL--MPC significantly outperformed standalone MPC with the same prediction horizon, indicating that the framework can mitigate the myopic behavior of short-horizon MPC.
When training a generalist policy on a broader set of weather disturbances and evaluating controller performance on weather trajectories held out during training, trajectory-selection RL--MPC improved over both standalone RL (by $54\%$) and MPC (by $80\%$).
Trajectory-selection RL--MPC achieved similar revenue from crop growth as the generalist policy, while reducing operating costs and constraint violations.

Future work should evaluate this framework over full tomato production cycles, using crop models that represent fruit-development stages, to assess whether this method can effectively balance current control decisions against longer-delayed economic returns.
In addition, future studies should benchmark trajectory-selection RL--MPC against control methods with longer prediction horizons that capture slow crop yield dynamics, such as the timescale decomposition approach, to quantify the trade-off between computational demand and closed-loop economic performance.
Finally, practical greenhouse experiments should assess the robustness of this control method under real-world uncertain conditions, such as model mismatch, state estimation errors, and unseen weather and market-price disturbances.

The source code required to reproduce the results is publicly available at \url{https://github.com/BartvLaatum/GL-Gym-MPC}.
Trained models and data for visualizations are available upon request.

\printcredits
\appendix

\begin{longtable}{p{0.52\linewidth}p{0.26\linewidth}p{0.16\linewidth}}
\caption{GreenLight model state variables with symbols and units.}\label{tab:gl_states}\\
\toprule
\textbf{State (description)} & \textbf{Symbol} & \textbf{Units} \\
\midrule
\endfirsthead

\toprule
\textbf{State (description)} & \textbf{Symbol} & \textbf{Units} \\
\midrule
\endhead

\midrule
\multicolumn{3}{r}{\emph{Continued on next page}}\\
\endfoot

\bottomrule
\endlastfoot

\multicolumn{3}{l}{\textbf{CO$_2$ variables}}\\
\midrule
CO$_2$ concentration of the main compartment            & $x_{\mathrm{CO_2},\mathrm{air}}$ & mg m$^{-3}$ s$^{-1}$ \\
CO$_2$ concentration of the top compartment             & $x_{\mathrm{CO_2},\mathrm{top}}$ & mg m$^{-3}$ s$^{-1}$ \\[2pt]
\midrule

\multicolumn{3}{l}{\textbf{Indoor temperature variables}}\\
\midrule
Greenhouse air temperature                              & $x_{\mathrm{T, air}}$            & $^\circ$C \\
Air temperature above the screen (top air)              & $x_{\mathrm{T,top}}$            & $^\circ$C \\
Canopy temperature                                      & $x_{\mathrm{T,can}}$            & $^\circ$C \\
Internal cover temperature                              & $x_{\mathrm{T,cov,in}}$       & $^\circ$C \\
External cover temperature                              & $x_{\mathrm{T,cov,out}}$      & $^\circ$C \\
Thermal screen temperature                              & $x_{\mathrm{T,thScr}}$            & $^\circ$C \\
Greenhouse floor temperature                            & $x_{\mathrm{T,flr}}$          & $^\circ$C \\
Pipe temperature                                        & $x_{\mathrm{T,pipe}}$           & $^\circ$C \\
First soil layer temperature                            & $x_{\mathrm{T,soil_1}}$       & $^\circ$C \\
Second soil layer temperature                           & $x_{\mathrm{T,soil_2}}$       & $^\circ$C \\
Third soil layer temperature                            & $x_{\mathrm{T,soil_3}}$       & $^\circ$C \\
Fourth soil layer temperature                           & $x_{\mathrm{T,soil_4}}$       & $^\circ$C \\
Fifth soil layer temperature                            & $x_{\mathrm{T,soil_5}}$       & $^\circ$C \\
Top lamp temperature                                    & $x_{\mathrm{T,lamp,top}}$     & $^\circ$C \\
Inter-lamp temperature                                  & $x_{\mathrm{T,lamp,inter}}$   & $^\circ$C \\
Grow pipe temperature                                   & $x_{\mathrm{T,growPipe}}$       & $^\circ$C \\
Blackout screen temperature                             & $x_{\mathrm{T,blScr}}$            & $^\circ$C \\[2pt]
\midrule

\multicolumn{3}{l}{\textbf{Vapour pressure variables}}\\
\midrule
Vapour pressure in the main compartment                 & $x_{\mathrm{VP},\mathrm{air}}$ & Pa \\
Vapour pressure in the top compartment                  & $x_{\mathrm{VP},\mathrm{top}}$ & Pa \\[2pt]
\midrule

\multicolumn{3}{l}{\textbf{Crop variables}}\\
\midrule
Average crop canopy temperature over the last 24 hours  & $x_{\mathrm{T},{\mathrm{can},\mathrm{24h}}}$ & $^\circ$C \\
Carbohydrates in the crop's buffer                      & $x_{\mathrm{buf}}$       & mg CH$_2$O m$^{-2}$ \\
Carbohydrates in the leaves                             & $x_{\mathrm{leaf}}$      & mg CH$_2$O m$^{-2}$ \\
Carbohydrates in the stem                               & $x_{\mathrm{stem}}$      & mg CH$_2$O m$^{-2}$ \\
Carbohydrates in the fruits                             & $x_{\mathrm{frt}}$     & mg CH$_2$O m$^{-2}$ \\
Crop temperature sum                                    & $x_{\mathrm{T},\mathrm{sum}}$       & $^\circ$C\,d$^{-1}$ \\
\end{longtable}

\bibliographystyle{cas-model2-names}

\bibliography{cas-refs}

@phdthesis{dekoningDevelopmentDryMatter1994,
    author = {De Koning, A.N.M.},
    title = {Development and Dry Matter Distribution in Glasshouse Tomato : A Quantitative Approach},
    school = {Agricultural University},
    year = 1994,
    doi = {10.18174/205947},
}

@article{adamsEffectTemperatureGrowth2001,
  title = {Effect of {{Temperature}} on the {{Growth}} and {{Development}} of {{Tomato Fruits}}},
  author = {Adams, S and Cave, C. R. J. and Cockshul, K. E.},
  year = {2001},
  journal = {Annals of Botany},
  volume = {88},
  number = {5},
  pages = {869--877},
  issn = {03057364},
  doi = {10.1006/anbo.2001.1524},
  urldate = {2026-05-31},
  langid = {english},
}

@article{vanstratenOptimalGreenhouseCultivation2010,
  title = {Optimal {{Greenhouse Cultivation Control}}: {{Survey}} and {{Perspectives}}},
  shorttitle = {Optimal {{Greenhouse Cultivation Control}}},
  author = {Van Straten, G. and Van Henten, E. J.},
  year` = {2010},
  journal = {IFAC Proceedings Volumes},
  shortjournal = {IFAC Proceedings Volumes},
  volume = {43},
  number = {26},
  pages = {18--33},
  issn = {14746670},
  doi = {10.3182/20101206-3-JP-3009.00004},
}

@phdthesis{tapEconomicsbasedOptimalControl2000,
      title = {Economics-Based Optimal Control of Greenhouse Tomato Crop Production},
  author = {Tap, F.},
  year = 2000,
  school = {Wageningen University and Research},
  doi = {10.18174/195235},
}

@article{kuijpersFruitDevelopmentModelling2021,
  title = {Fruit Development Modelling and Performance Analysis of Automatic Greenhouse Control},
  author = {Kuijpers, Wouter J.P. and Antunes, Duarte J. and Hemming, Silke and Van Henten, Eldert J. and Van De Molengraft, Marinus J.G.},
  year = {2021},
  journal = {Biosystems Engineering},
  shortjournal = {Biosystems Engineering},
  volume = {208},
  pages = {300--318},
  issn = {15375110},
  doi = {10.1016/j.biosystemseng.2021.06.002},
}

@article{VANHENTEN2009timescale,
  title = {Time-Scale Decomposition of an Optimal Control Problem in Greenhouse Climate Management},
  author = {Van Henten, E.J. and Bontsema, J.},
  year = {2009},
  journal = {Control Engineering Practice},
  volume = {17},
  number = {1},
  pages = {88--96},
  issn = {0967-0661},
  doi = {10.1016/j.conengprac.2008.05.008},
}

@article{kearns2002near,
  title={Near-optimal reinforcement learning in polynomial time},
  author={Kearns, Michael and Singh, Satinder},
  journal={Machine learning},
  volume={49},
  number={2},
  pages={209--232},
  year={2002},
  doi={https://doi.org/10.1023/A:1017984413808},
  publisher={Springer}
}

@inproceedings{serban2005cvodes,
  title={CVODES: the sensitivity-enabled ODE solver in SUNDIALS},
  author={Serban, Radu and Hindmarsh, Alan C},
  booktitle={International Design Engineering Technical Conferences and Computers and Information in Engineering Conference},
  volume={47438},
  pages={257--269},
  year={2005},
  doi={https://doi.org/10.1115/DETC2005-85597}
}

@article{anderssonCasADiSoftwareFramework2019a,
  title = {{{CasADi}}: A Software Framework for Nonlinear Optimization and Optimal Control},
  shorttitle = {{{CasADi}}},
  author = {Andersson, Joel A. E. and Gillis, Joris and Horn, Greg and Rawlings, James B. and Diehl, Moritz},
  year = 2019,
  month = mar,
  journal = {Mathematical Programming Computation},
  volume = {11},
  number = {1},
  pages = {1--36},
  issn = {1867-2957},
  doi = {10.1007/s12532-018-0139-4},
  urldate = {2025-07-02},
  langid = {english},
}

@article{wachterImplementationInteriorpointFilter2006,
  title = {On the Implementation of an Interior-Point Filter Line-Search Algorithm for Large-Scale Nonlinear Programming},
  author = {W{\"a}chter, Andreas and Biegler, Lorenz T.},
  year = 2006,
  month = mar,
  journal = {Mathematical Programming},
  volume = {106},
  number = {1},
  pages = {25--57},
  issn = {1436-4646},
  doi = {10.1007/s10107-004-0559-y},
  urldate = {2025-07-02},
  langid = {english},
}

@article{raffinStableBaselines3ReliableReinforcement2021,
  title = {Stable-{{Baselines3}}: {{Reliable Reinforcement Learning Implementations}}},
  shorttitle = {Stable-{{Baselines3}}},
  author = {Raffin, Antonin and Hill, Ashley and Gleave, Adam and Kanervisto, Anssi and Ernestus, Maximilian and Dormann, Noah},
  year = 2021,
  journal = {Journal of Machine Learning Research},
  volume = {22},
  number = {268},
  pages = {1--8},
  issn = {1533-7928},
  urldate = {2025-07-02},
}

@misc{schulmanProximalPolicyOptimization2017,
  title = {Proximal {{Policy Optimization Algorithms}}},
  author = {Schulman, John and Wolski, Filip and Dhariwal, Prafulla and Radford, Alec and Klimov, Oleg},
  year = 2017,
  month = aug,
  number = {arXiv:1707.06347},
  eprint = {1707.06347},
  primaryclass = {cs},
  publisher = {arXiv},
  doi = {10.48550/arXiv.1707.06347},
  urldate = {2025-11-05},
  archiveprefix = {arXiv},
  langid = {english},
}

@article{vanthoorMethodologyModelbasedGreenhouse2011,
  title = {A Methodology for Model-Based Greenhouse Design: {{Part}} 2, Description and Validation of a Tomato Yield Model},
  shorttitle = {A Methodology for Model-Based Greenhouse Design},
  author = {Vanthoor, B.H.E. and De Visser, P.H.B. and Stanghellini, C. and Van Henten, E.J.},
  year = 2011,
  month = dec,
  journal = {Biosystems Engineering},
  volume = {110},
  number = {4},
  pages = {378--395},
  issn = {15375110},
  doi = {10.1016/j.biosystemseng.2011.08.005},
  urldate = {2025-02-26},
  copyright = {https://www.elsevier.com/tdm/userlicense/1.0/},
  langid = {english},
}

@misc{mediratta2024generalizationgapofflinereinforcement,
      title={The Generalization Gap in Offline Reinforcement Learning}, 
      author={Ishita Mediratta and Qingfei You and Minqi Jiang and Roberta Raileanu},
      year={2024},
      eprint={2312.05742},
      archivePrefix={arXiv},
      primaryClass={cs.LG},
      url={https://arxiv.org/abs/2312.05742}, 
}

@article{kirkSurveyZeroshotGeneralisation2023,
  title = {A {{Survey}} of {{Zero-shot Generalisation}} in {{Deep Reinforcement Learning}}},
  author = {Kirk, Robert and Zhang, Amy and Grefenstette, Edward and Rockt{\"a}schel, Tim},
  year = 2023,
  month = jan,
  journal = {Journal of Artificial Intelligence Research},
  volume = {76},
  pages = {201--264},
  issn = {1076-9757},
  doi = {10.1613/jair.1.14174},
  urldate = {2025-02-13},
  langid = {english},
}

@misc{zhangDissectionOverfittingGeneralization2018,
  title = {A {{Dissection}} of {{Overfitting}} and {{Generalization}} in {{Continuous Reinforcement Learning}}},
  author = {Zhang, Amy and Ballas, Nicolas and Pineau, Joelle},
  year = 2018,
  month = jun,
  number = {arXiv:1806.07937},
  eprint = {1806.07937},
  primaryclass = {cs},
  publisher = {arXiv},
  doi = {10.48550/arXiv.1806.07937},
  urldate = {2025-02-13},
  archiveprefix = {arXiv},
  langid = {english},
}

@misc{zhangStudyOverfittingDeep2018,
  title = {A {{Study}} on {{Overfitting}} in {{Deep Reinforcement Learning}}},
  author = {Zhang, Chiyuan and Vinyals, Oriol and Munos, Remi and Bengio, Samy},
  year = 2018,
  month = apr,
  number = {arXiv:1804.06893},
  eprint = {1804.06893},
  primaryclass = {cs},
  publisher = {arXiv},
  doi = {10.48550/arXiv.1804.06893},
  urldate = {2026-03-31},
  archiveprefix = {arXiv},
  langid = {english},
}

@article{katzinGreenLightOpenSource2020,
  title = {{{GreenLight}} -- {{An}} Open Source Model for Greenhouses with Supplemental Lighting: {{Evaluation}} of Heat Requirements under {{LED}} and {{HPS}} Lamps},
  shorttitle = {{{GreenLight}} -- {{An}} Open Source Model for Greenhouses with Supplemental Lighting},
  author = {Katzin, David and Van Mourik, Simon and Kempkes, Frank and Van Henten, Eldert J.},
  year = 2020,
  month = jun,
  journal = {Biosystems Engineering},
  volume = {194},
  pages = {61--81},
  issn = {15375110},
  doi = {10.1016/j.biosystemseng.2020.03.010},
  urldate = {2025-01-31},
  langid = {english},
}

@inproceedings{MSAAD2025449,
  title = {{{RL-guided MPC}} for Autonomous Greenhouse Control},
  booktitle = {{{IFAC-PapersOnLine}}},
  author = {Msaad, Salim and Harraway, Murray and McAllister, Robert D.},
  year = 2025,
  volume = {59},
  pages = {449--454},
  doi = {10.1016/j.ifacol.2025.11.829}
}

@article{vanlaatumStochasticModelPredictive2026,
  title = {Stochastic Model Predictive Control with Reinforcement Learning for Greenhouse Production Systems under Parametric Uncertainty},
  author = {Van Laatum, Bart and Msaad, Salim and Van Henten, Eldert J. and Mcallister, Robert D. and Boersma, Sjoerd},
  year = 2026,
  month = apr,
  journal = {Control Engineering Practice},
  volume = {169},
  pages = {106787},
  issn = {09670661},
  doi = {10.1016/j.conengprac.2026.106787},
  urldate = {2026-02-03},
  langid = {english},
}

@article{reiter2025ac4mpc,
  title={AC4MPC: Actor-Critic Reinforcement Learning for Guiding Model Predictive Control},
  author={Reiter, Rudolf and Ghezzi, Andrea and Baumg{\"a}rtner, Katrin and Hoffmann, Jasper and McAllister, Robert D and Diehl, Moritz},
  journal={IEEE Transactions on Control Systems Technology},
  year={2025},
  publisher={IEEE}
}

@misc{ghezzi2026rolloutoptimizeonestepnewton,
      title={Rollout Then Optimize: A One-Step Newton Refinement of Learned Policies for Nonlinear Model Predictive Control}, 
      author={Andrea Ghezzi and Rudolf Reiter and Katrin Baumgärtner and Alberto Bemporad and Moritz Diehl},
      year={2026},
      eprint={2504.02710},
      archivePrefix={arXiv},
      primaryClass={math.OC},
      url={https://arxiv.org/abs/2504.02710}, 
}

@inproceedings{vanlaatum2025greenlight,
  title={GreenLight-Gym: Reinforcement learning benchmark environment for control of greenhouse production systems},
  author={Van Laatum, Bart and Van Henten, Eldert J and Boersma, Sjoerd},
  booktitle={8th IFAC Conference on Sensing, Control and Automation Technologies for Agriculture, AGRICONTROL 2025},
  pages={437--442},
  year={2025},
  doi={https://doi.org/10.1016/j.ifacol.2025.11.827},
  organization={IFAC}
}

@book{suttonReinforcementLearningIntroduction2018a,
  title = {Reinforcement Learning: An Introduction},
  shorttitle = {Reinforcement Learning},
  author = {Sutton, Richard S. and Barto, Andrew G.},
  year = 2018,
  series = {Adaptive Computation and Machine Learning Series},
  edition = {Second edition},
  publisher = {The MIT Press},
  address = {Cambridge, Massachusetts},
  isbn = {978-0-262-03924-6},
  langid = {english},
  lccn = {Q325.6 .R45 2018},
}

@article{blascoModelbasedPredictiveControl2007,
  title = {Model-Based Predictive Control of Greenhouse Climate for Reducing Energy and Water Consumption},
  author = {Blasco, X. and Mart{\'i}nez, M. and Herrero, J. M. and Ramos, C. and Sanchis, J.},
  year = 2007,
  month = jan,
  journal = {Computers and Electronics in Agriculture},
  volume = {55},
  number = {1},
  pages = {49--70},
  issn = {0168-1699},
  doi = {10.1016/j.compag.2006.12.001},
  urldate = {2025-07-02}
}

@article{boersmaRobustSamplebasedModel2022,
  title = {Robust Sample-Based Model Predictive Control of a Greenhouse System with Parametric Uncertainty},
  author = {Boersma, Sjoerd and Sun, Congcong and Van Mourik, Simon},
  year = 2022,
  journal = {IFAC-PapersOnLine},
  volume = {55},
  number = {32},
  pages = {177--182},
  issn = {24058963},
  doi = {10.1016/j.ifacol.2022.11.135},
  urldate = {2024-12-02},
  langid = {english}
}

@inproceedings{finkComparisonDynamicTomato2023,
  title = {Comparison of {{Dynamic Tomato Growth Models}} for {{Optimal Control}} in {{Greenhouses}}},
  booktitle = {2023 {{IEEE International Conference}} on {{Agrosystem Engineering}}, {{Technology}} \& {{Applications}} ({{AGRETA}})},
  author = {Fink, Michael and Daniels, Annalena and Qian, Cheng and Vel{\'a}squez, V{\'i}ctor Mart{\'i}nez and Salotra, Sahil and Wollherr, Dirk},
  year = 2023,
  month = sep,
  pages = {28--33},
  publisher = {IEEE},
  address = {Shah Alam, Malaysia},
  doi = {10.1109/AGRETA57740.2023.10262422},
  urldate = {2026-03-30},
  copyright = {https://doi.org/10.15223/policy-029},
  isbn = {979-8-3503-4733-3},
  langid = {english}
}

@article{gruberNonlinearMPCBased2011,
  title = {Nonlinear {{MPC}} Based on a {{Volterra}} Series Model for Greenhouse Temperature Control Using Natural Ventilation},
  author = {Gruber, J. K. and Guzm{\'a}n, J. L. and Rodr{\'i}guez, F. and Bordons, C. and Berenguel, M. and S{\'a}nchez, J. A.},
  year = 2011,
  month = apr,
  journal = {Control Engineering Practice},
  volume = {19},
  number = {4},
  pages = {354--366},
  issn = {0967-0661},
  doi = {10.1016/j.conengprac.2010.12.004},
  urldate = {2025-07-02}
}

@inproceedings{itoGreenhouseTemperatureControl2012,
  title = {Greenhouse Temperature Control with Wooden Pellet Heater via Model Predictive Control Approach},
  booktitle = {2012 20th {{Mediterranean Conference}} on {{Control}} \& {{Automation}} ({{MED}})},
  author = {Ito, Kazuhisa},
  year = 2012,
  month = jul,
  pages = {1542--1547},
  doi = {10.1109/MED.2012.6265858},
  urldate = {2025-07-02}
}

@article{kuijpersLightingSystemsStrategies2021,
  title = {Lighting Systems and Strategies Compared in an Optimally Controlled Greenhouse},
  author = {Kuijpers, Wouter J.P. and Katzin, David and Van Mourik, Simon and Antunes, Duarte J. and Hemming, Silke and Van De Molengraft, Marinus J.G.},
  year = 2021,
  month = feb,
  journal = {Biosystems Engineering},
  volume = {202},
  pages = {195--216},
  issn = {15375110},
  doi = {10.1016/j.biosystemseng.2020.12.006},
  urldate = {2024-12-03},
  langid = {english}
}

@article{mallickReinforcementLearningbasedModel2025,
  title = {Reinforcement Learning-Based Model Predictive Control for Greenhouse Climate Control},
  author = {Mallick, Samuel and Airaldi, Filippo and Dabiri, Azita and Sun, Congcong and De Schutter, Bart},
  year = 2025,
  month = mar,
  journal = {Smart Agricultural Technology},
  volume = {10},
  pages = {100751},
  issn = {27723755},
  doi = {10.1016/j.atech.2024.100751},
  urldate = {2025-06-30},
  langid = {english}
}

@article{montoyaHybridcontrolledApproachMaintaining2016,
  title = {A Hybrid-Controlled Approach for Maintaining Nocturnal Greenhouse Temperature: {{Simulation}} Study},
  shorttitle = {A Hybrid-Controlled Approach for Maintaining Nocturnal Greenhouse Temperature},
  author = {Montoya, A. P. and Guzm{\'a}n, J. L. and Rodr{\'i}guez, F. and {S{\'a}nchez-Molina}, J. A.},
  year = 2016,
  month = apr,
  journal = {Computers and Electronics in Agriculture},
  volume = {123},
  pages = {116--124},
  issn = {0168-1699},
  doi = {10.1016/j.compag.2016.02.014},
  urldate = {2025-07-02}
}

@article{morcegoReinforcementLearningModel2023,
  title = {Reinforcement {{Learning}} versus {{Model Predictive Control}} on Greenhouse Climate Control},
  author = {Morcego, Bernardo and Yin, Wenjie and Boersma, Sjoerd and Van Henten, Eldert and Puig, Vicen{\c c} and Sun, Congcong},
  year = 2023,
  month = dec,
  journal = {Computers and Electronics in Agriculture},
  volume = {215},
  pages = {108372},
  issn = {01681699},
  doi = {10.1016/j.compag.2023.108372},
  urldate = {2025-02-13},
  langid = {english}
}

@article{ramirez-ariasMultiobjectiveHierarchicalControl2012,
  title = {Multiobjective Hierarchical Control Architecture for Greenhouse Crop Growth},
  author = {{Ram{\'i}rez-Arias}, A. and Rodr{\'i}guez, F. and Guzm{\'a}n, J.L. and Berenguel, M.},
  year = 2012,
  month = mar,
  journal = {Automatica},
  volume = {48},
  number = {3},
  pages = {490--498},
  issn = {00051098},
  doi = {10.1016/j.automatica.2012.01.002},
  urldate = {2026-03-30},
  copyright = {https://www.elsevier.com/tdm/userlicense/1.0/},
  langid = {english}
}

@article{svensenChanceconstrainedStochasticMPC2024,
  title = {Chance-Constrained Stochastic {{MPC}} of Greenhouse Production Systems with Parametric Uncertainty},
  author = {Svensen, Jan Lorenz and Cheng, Xiaodong and Boersma, Sjoerd and Sun, Congcong},
  year = 2024,
  month = feb,
  journal = {Computers and Electronics in Agriculture},
  volume = {217},
  pages = {108578},
  issn = {01681699},
  doi = {10.1016/j.compag.2023.108578},
  urldate = {2024-12-02},
  langid = {english}
}

@phdthesis{vanhentenGreenhouseClimateManagement1994,
  title = {{Greenhouse climate management : an optimal control approach}},
  shorttitle = {{Greenhouse climate management}},
  author = {Van Henten, E.J.},
  year = 1994,
  doi = {10.18174/205106},
  urldate = {2025-01-29},
  langid = {dutch},
  school = {Agricultural University},
}

@book{vanstratenOptimalControlGreenhouse2010,
  title = {Optimal Control of Greenhouse Cultivation},
  author = {Van Straten, Gerrit and Van Willigenburg, Gerard and Van Henten, E.J. and Van Ooteghem, Rachel},
  year = {2010},
  month = nov,
  edition = {0},
  publisher = {CRC Press},
  doi = {10.1201/b10321},
  urldate = {2025-02-13},
  isbn = {978-0-429-13731-0},
  langid = {english}
}

@article{xuAdaptiveTwoTimescale2018,
  title = {Adaptive Two Time-Scale Receding Horizon Optimal Control for Greenhouse Lettuce Cultivation},
  author = {Xu, Dan and Du, Shangfeng and Van Willigenburg, Gerard},
  year = 2018,
  month = mar,
  journal = {Computers and Electronics in Agriculture},
  volume = {146},
  pages = {93--103},
  issn = {01681699},
  doi = {10.1016/j.compag.2018.02.001},
  urldate = {2024-12-02},
  langid = {english}
}

@article{xuRulebasedYearroundModel2025,
  title = {Rule-Based Year-Round Model Predictive Control of Greenhouse Tomato Cultivation: {{A}} Simulation Study},
  shorttitle = {Rule-Based Year-Round Model Predictive Control of Greenhouse Tomato Cultivation},
  author = {Xu, Dan and Xu, Lei and Wang, Shusheng and Wang, Mingqin and Ma, Juncheng and Shi, Chen},
  year = 2025,
  month = sep,
  journal = {Information Processing in Agriculture},
  volume = {12},
  number = {3},
  pages = {344--357},
  issn = {22143173},
  doi = {10.1016/j.inpa.2024.11.001},
  urldate = {2026-03-30},
  langid = {english}
}

@misc{knmi_uurgeg_schiphol_2011_2020,
  author       = {{KNMI}},
  title        = {Uurgegevens van het weer in Nederland: station 240 (Schiphol), 2011--2020},
  year         = {2020},
  howpublished = {ZIP archive (ASCII)},
  url          = {https://cdn.knmi.nl/knmi/map/page/klimatologie/gegevens/uurgegevens/uurgeg_240_2011-2020.zip},
  copyright     = {CC BY 4.0},
  note         = {Accessed 2025. Open data provided under CC BY 4.0;}
}

@misc{agc_2023,
  author = {Isabella Righini and Bart van Marrewijk and Georgios Ntakos and Pinglin Zhang and Guido Jansen and S.C. Maree and Monique Bijlaard  and H. F. de Zwart and Silke Hemming},
  title = {4th Autonomous Greenhouse Challenge: Pre-trial Dwarf Tomato Measurements and
  Images},
  doi = {10.4121/e1ee9de9-6ce9-4502-a37c-34b5b1372bed.v1},
  publisher = {4TU.ResearchData},
  year = {2024},
  copyright = {CC BY 4.0},
  note         = {Accessed 2025. Open data provided under CC BY 4.0;}
}

@article{petropoulou2023lettuce,
AUTHOR = {Petropoulou, Anna Selini and van Marrewijk, Bart and de Zwart, Feije and Elings, Anne and Bijlaard, Monique and van Daalen, Tim and Jansen, Guido and Hemming, Silke},
TITLE = {Lettuce Production in Intelligent Greenhouses—3D Imaging and Computer Vision for Plant Spacing Decisions},
JOURNAL = {Sensors},
VOLUME = {23},
YEAR = {2023},
NUMBER = {6},
ARTICLE-NUMBER = {2929},
PubMedID = {36991638},
ISSN = {1424-8220},
DOI = {10.3390/s23062929}
}

@article{maree2025autonomous,
AUTHOR = {Maree, Stef C. and Zhang, Pinglin and van Marrewijk, Bart M. and de Zwart, Feije and Bijlaard, Monique and Hemming, Silke},
TITLE = {Autonomous Greenhouse Cultivation of Dwarf Tomato: Performance Evaluation of Intelligent Algorithms for Multiple-Sensor Feedback},
JOURNAL = {Sensors},
VOLUME = {25},
YEAR = {2025},
NUMBER = {14},
ARTICLE-NUMBER = {4321},
PubMedID = {40732449},
ISSN = {1424-8220},
DOI = {10.3390/s25144321}
}

@article{ariesen-verschuurDigitalTwinsGreenhouse2022,
  title = {Digital {{Twins}} in Greenhouse Horticulture: {{A}} Review},
  shorttitle = {Digital {{Twins}} in Greenhouse Horticulture},
  author = {{Ariesen-Verschuur}, Natasja and Verdouw, Cor and Tekinerdogan, Bedir},
  year = 2022,
  month = aug,
  journal = {Computers and Electronics in Agriculture},
  volume = {199},
  pages = {107183},
  issn = {01681699},
  doi = {10.1016/j.compag.2022.107183},
  urldate = {2026-03-29},
  langid = {english}
}

@article{bisbisPotentialImpactsClimate2018,
  title = {Potential Impacts of Climate Change on Vegetable Production and Product Quality -- {{A}} Review},
  author = {Bisbis, Mehdi Benyoussef and Gruda, Nazim and Blanke, Michael},
  year = 2018,
  month = jan,
  journal = {Journal of Cleaner Production},
  volume = {170},
  pages = {1602--1620},
  issn = {09596526},
  doi = {10.1016/j.jclepro.2017.09.224},
  urldate = {2026-03-29},
  langid = {english}
}

@article{fanourakisClimateChangeImpacts2025,
  title = {Climate {{Change Impacts}} on {{Greenhouse Horticulture}} in the {{Mediterranean Basin}}: {{Challenges}} and {{Adaptation Strategies}}},
  shorttitle = {Climate {{Change Impacts}} on {{Greenhouse Horticulture}} in the {{Mediterranean Basin}}},
  author = {Fanourakis, Dimitrios and Tsaniklidis, Georgios and Makraki, Theodora and Nikoloudakis, Nikolaos and Bartzanas, Thomas and Sabatino, Leo and Fatnassi, Hicham and Ntatsi, Georgia},
  year = 2025,
  month = nov,
  journal = {Plants},
  volume = {14},
  number = {21},
  pages = {3390},
  issn = {2223-7747},
  doi = {10.3390/plants14213390},
  urldate = {2026-03-29},
  langid = {english}
}

@article{muldersExtremeDroughtRainfall2024,
  title = {Extreme Drought and Rainfall Had a Large Impact on Potato Production in the {{Netherlands}} between 2015 and 2020},
  author = {Mulders, Puck J. A. M. and Van Den Heuvel, Edwin R. and Van De Molengraft, M. J. G. and Heemels, W. P. M. H. and Reidsma, Pytrik},
  year = 2024,
  month = sep,
  journal = {Communications Earth \& Environment},
  volume = {5},
  number = {1},
  pages = {496},
  issn = {2662-4435},
  doi = {10.1038/s43247-024-01658-3},
  urldate = {2026-03-29},
  langid = {english}
}

@book{faoGoodAgriculturalPratices2013,
  title = {Good Agricultural Pratices for Greenhouse Vegetable Crops: Principles for Mediterranean Climate Areas},
  shorttitle = {Good Agricultural Pratices for Greenhouse Vegetable Crops},
  author={Baudoin, Wilfried and Nono-Womdim, Remi and Lutaladio, NeBambi and Hodder, Alison and Castilla, Nicol{\'a}s and Leonardi, Cherubino and De Pascale, Stefania and Qaryouti, Muien and Duffy, R},
  year = {2013},
  number = {217},
  publisher = {FAO},
  address = {Rome},
  isbn = {978-92-5-107649-1},
  langid = {english},
  lccn = {635.982 3},
}

@article{silvaAssessingImpactGlobal2017,
  title = {Assessing the Impact of Global Warming on Worldwide Open Field Tomato Cultivation through {{CSIRO-Mk3}}{$\cdot$}0 Global Climate Model},
  author = {Silva, R. S. and Kumar, L. and Shabani, F. and Pican{\c c}o, M. C.},
  year = 2017,
  month = apr,
  journal = {The Journal of Agricultural Science},
  volume = {155},
  number = {3},
  pages = {407--420},
  issn = {0021-8596, 1469-5146},
  doi = {10.1017/S0021859616000654},
  urldate = {2026-03-29},
  copyright = {https://www.cambridge.org/core/terms},
  langid = {english}
}



\end{document}